\documentclass[final,hidelinks,onefignum,onetabnum]{siamart250211}

\usepackage{booktabs}        
\usepackage{mathtools}       
\usepackage{microtype}       
\usepackage{graphicx}        
\usepackage{enumitem}        
\usepackage{url}
\usepackage{xcolor}
\usepackage{amsfonts}
\usepackage{algorithm}
\usepackage{algorithmic}
\usepackage{amssymb}

\providecommand{\citet}[1]{\cite{#1}}
\providecommand{\citep}[1]{\cite{#1}}
 
\providecommand{\citeauthor}[1]{\textsc{#1}}
\providecommand{\citeyear}[1]{\cite{#1}}

\providecommand{\tr}{\mathrm{tr}}

\graphicspath{{./}}

\newtheorem{remark}{Remark}

\theoremstyle{plain}

\DeclareMathOperator{\spec}{Spec}

\headers{BOLT: Block-Orthonormal Lanczos for trace estimation}{K. Yeon, P. Ghosal, and M. Anitescu}

\title{BOLT: Block-Orthonormal Lanczos for Trace Estimation of Matrix Functions}

\author{
  Kingsley Yeon\thanks{Department of Statistics, University of Chicago
  (\email{yeon@uchicago.edu}).}
  \and
  Promit Ghosal\thanks{(\email{promit@uchicago.edu}).}
  \and
  Mihai Anitescu\thanks{Mathematics and Computer Science, Argonne National Laboratory
  (\email{anitescu@mcs.anl.gov}).}
}

\begin{document}
\maketitle

\begin{abstract}
Efficient matrix trace estimation is essential for scalable computation of log-determinants,
matrix norms, and distributional divergences. In many large-scale applications, the matrices
involved are too large to store or access in full, making even a single matrix--vector (mat--vec)
product infeasible. Instead, one often has access only to small subblocks of the matrix or
localized matrix--vector products on restricted index sets. Hutch++~\cite{meyer2021hutch++}
achieves an optimal convergence rate but relies on randomized SVD and assumes full mat--vec
access, making it difficult to apply in these constrained settings. We propose the
Block-Orthonormal Stochastic Lanczos Quadrature (BOLT), which outperforms Hutch++
in near flat-spectrum settings while maintaining a simpler implementation based on
orthonormal block probes and Lanczos iterations. BOLT builds on the Stochastic Lanczos
Quadrature (SLQ) framework, which combines random probing with Krylov subspace methods
to efficiently approximate traces of matrix functions. To address memory limitations and
partial access constraints, we introduce Subblock SLQ, a variant of BOLT that operates only
on small principal submatrices. In a large-scale example, we recover the trace
of a \(10^6 \times 10^6\) matrix of the form \(A = B^T B\) with
\(B \in \mathbb{R}^{2048 \times 10^6}\) by sampling diagonal subblocks of size 64, achieving
a relative error of \(3.78\times 10^{-5}\) while observing less than 10\% of the matrix
entries. As a result, this framework yields a proxy KL divergence estimator and an efficient
method for computing the Wasserstein-2 distance between Gaussians, both compatible with
low-memory and partial-access regimes. We provide theoretical guarantees and demonstrate
strong empirical performance across a range of high-dimensional settings.
\end{abstract}



\section{Introduction}
Trace estimation is central to scientific computing, enabling approximation of log-determinants~\cite{boutsidis2017randomized}, matrix norms~\cite{han2017approximating}, distributional divergences, spectral densities~\cite{lin2016approximating}, the Estrada index~\cite{wang2021public}, and more~\cite{ubaru2017applications}. In settings such as covariance estimation, kernel learning, and Gaussian processes, efficient trace estimation underlies key computational pipelines~\cite{cai2015law, dong2017scalable, wenger2022preconditioning}.

Recent randomized methods have made trace estimation highly efficient. Hutchinson’s
estimator uses random sign probes and requires \(O(1/\epsilon^2)\) matrix–vector products
for additive error \(\epsilon\) with high probability~\cite{avron2011randomized}, while Hutch++
reduces this to \(O(1/\epsilon)\) via low-rank sketching~\cite{meyer2021hutch++}. Non-adaptive
sketching~\cite{jiang2021optimal} improves memory and parallelism, establishing this rate as
optimal. SLQ achieves fast convergence for analytic functions~\cite{ubaru2017fast}, and
ContHutch++ generalizes to implicit operators~\cite{zvonek2025conthutchplusplus}. However,
many of these methods still require full matrix–vector products and do not fully exploit
sparsity or localized access patterns.

Our work is motivated by estimating the KL divergence between mean-zero Gaussians in high-dimensional settings, where covariance matrices are large, possibly singular, and often inaccessible in full. These challenges arise in precision matrix estimation for Kalman filtering in weather data assimilation, where memory limits prohibit full matrix–vector products and the sample covariance is rank-deficient due to limited observations~\cite{fan2011high, lounici2014high, furrer2007estimation}. In such regimes, the KL divergence is not well-defined, motivating subblock-based approaches.

To address this, we propose Subblock SLQ—a memory-aware variant of BOLT that achieves fast convergence using only small principal submatrices. This method supports trace estimation without full matrix access at a time and remains effective even when the covariance is singular. Within this framework, we define a proxy KL divergence that remains meaningful without full-rank assumptions, enabling KL-based comparisons in sample-starved regimes. When KL is not required, the same estimator also supports the Wasserstein-2 (\(W_2\)) distance, which avoids matrix inversion and applies to any positive semidefinite covariance. Both estimators rely entirely on localized computations, making them ideal for large-scale, memory-constrained applications.


\noindent
Our key contributions are as follows:

(1) We develop the \textit{Block-Orthonormal SLQ (BOLT)} algorithm, a block-based stochastic
Lanczos quadrature method that achieves accuracy with optimal FLOP complexity
(see Table~\ref{table:flop}). Unlike Hutch++~\cite{meyer2021hutch++}, BOLT avoids
randomized SVD, making it easier to implement and more robust in flat-spectrum
settings such as \(\mathrm{diag}(\mathrm{Unif}[1,2])\) (see Figure~\ref{fig:hutchppfail}),
where it outperforms Hutch++. BOLT leverages block probes well-suited to modern
hardware and BLAS-efficient matrix–matrix operations, and includes theoretical
guarantees for unbiasedness and variance control (Theorem~\ref{thm:blockSLQ}). Its
computational benefits are demonstrated in Figure~\ref{fig:error_vs_matvec_actual_N}.

(2) Building on BOLT, we propose the \textit{Subblock SLQ} method for estimating divergences using only small principal submatrices, enabling fast and memory-efficient trace computation when full mat–vecs are infeasible. This method naturally supports both a proxy KL divergence and the Wasserstein-2 (\(W_2\)) distance, with accuracy and efficiency governed by localized matrix structure. Theoretical guarantees are provided in Algorithm~\ref{algo:unified}, Theorem~\ref{thm:localization}, and Lemma~\ref{lem:unbiasedsubblockbound}.

(3) 
When the KL divergence is the target metric and the underlying covariance matrix is singular, we define a trace-based proxy for the KL divergence using Subblock SLQ, which is computable from small principal submatrices. This formulation remains well-defined even when the covariance matrix is singular. Lemma~\ref{thm:fullrank} shows that principal submatrices can retain full rank as long as their size does not exceed the rank of the full matrix, ensuring that the proxy KL divergence remains stable and tractable in regimes where the classical KL is undefined.


(4) 
We validate the proposed proxy KL divergence as a practical metric through various
experiments for assessing estimation error in low-sample regimes, particularly when the
covariance matrix is singular and the classical KL divergence is undefined. To this end,
we compare two strategies for estimating the Cholesky factor of the precision matrix:
one that enforces sparsity on a lower-triangular factor \(L\), and another that solves an
unconstrained least-squares problem. The latter approach is ill-posed and yields
substantial errors (Figure~\ref{fig:err_sparse}), which are accurately captured by our
proxy KL divergence. Furthermore, we demonstrate that SLQ-KL regularization
enhances generalization in low-rank learning tasks, including MNIST image
classification (Figure~\ref{fig:accuracy_sql}), by promoting isotropy—a property known
to improve generalization in both NLP and vision domains~\cite{arora2017simple,
mu2017all, cai2021isotropy}. Finally, we show that the same proxy KL framework is
useful for evaluating hierarchical off-diagonal low-rank (HODLR) approximations:
because there is no fast way to assess global errors such as \(\|A-\tilde A\|_F\) or verify
\(H^{-1}A \approx I_n\) without \(n^2\) access, subblock SLQ provides an efficient
alternative that quantifies the divergence between a kernel matrix and its HODLR
representation using only fast hierarchical mat–vecs
(Figures~\ref{fig:hodlr-peel-lvl1}–\ref{fig:hodlr-peel-lvl2}).

\subsection{Related Works and Connections}

\paragraph{Trace Estimation}
Stochastic trace estimation is well-studied~\cite{meyer2021hutch++, epperly2024xtrace, avron2011randomized}, with widely used methods such as Hutch++, XTrace, and XNysTrace offering effective variance reduction under favorable conditions. Hutch++, XTrace, and XNysTrace estimate the trace by sketching dominant eigenspaces, with XTrace and XNysTrace adding leave-one-out corrections using the exchangeability principle; these methods rely on full matrix–vector products and perform best under fast spectral decay, offering limited benefit for flat-spectrum matrices. XTrace and XNysTrace also incur higher computational overhead: their naïve implementations rely on repeated QR factorizations (XTrace) or SVDs (XNysTrace), costing $O(N m^3)$ and $O(m^4)$, respectively, where $m$ represents the number of matrix-vector products and $N$ denotes the size of the matrix ($N \times N$). Optimized variants reduce this to $O(m^3)$ by exploiting a rank-one update structure, but still demand greater infrastructure than simpler estimators. Table~\ref{table:flop} compares the implementations; see Appendix~\ref{appendix:trace_flops} for derivations and Table~\ref{table:flop-old} for details.

\begin{table}[h]
\centering
\caption{Flop counts to achieve \(\varepsilon\)-accuracy (\(m, b \sim \Theta(\varepsilon^{-1})\)) (\textbf{Optimized}).}
\label{table:flop}
\begin{tabular}{l c c c}
\toprule
Method      & \# Mat–vecs & Overhead (post–mat–vec)         & Total flops                              \\
\midrule
XTrace      & \(m\)  
            & \(O(N\,m^2 + m^3)\)             & \(O(N^2 m + N\,m^2 + m^3)\)               \\

XNysTrace   & \(m\)  
            & \(O(N\,m^2 + m^3)\)             & \(O(N^2 m + N\,m^2 + m^3)\)               \\

Hutch++     & \(m\)  
            & \(O(N\,m^2)\)                   & \(O(N^2 m + N\,m^2)\)                     \\

BOLT        & \(m\)  
            & \(O(N\,m^2)\)                   & \(O(N^2 m + N\,m^2)\)                     \\
\bottomrule
\end{tabular}
\end{table}

Fuentes et al.~\cite{fuentes2023estimating} empirically observed that block probe vectors yield better accuracy than scalar ones and noted rapid convergence due to self-averaging, though no variance bound was provided. Their analysis was limited to full-matrix settings with fixed block probes; they did not consider the role of orthogonalization, nor explore subblock strategies suitable for large-scale or structured matrices. None of these methods explicitly address memory usage or scalability under storage constraints, which are often critical in high-dimensional applications.

In a similar spirit as block SLQ, sliced-Wasserstein (SW) distances provides the framework to compute distance between two probability measure by projecting measures onto random lines:
\[
\mathrm{SW}_p^p(P,Q) = \int_{\mathbb{S}^{d-1}} W_p^p(P_\theta, Q_\theta)\,d\theta,
\]
where \(P_\theta\), \(Q_\theta\) are 1D marginals of \(P\), \(Q\) along \(\theta\)~\cite{nadjahi2021fast, kolouri2019generalized}. For example, if \(P = \mathcal{N}(0, uu^T)\), \(Q = \mathcal{N}(0, vv^T)\) are rank-one Gaussians, projections orthogonal to both \(u\) and \(v\) collapse to the origin, yielding zero SW distance despite distinct inputs. Empirically, SW distances are estimated via Monte Carlo integration:
\[
\widehat{\mathrm{SW}}_p^p(P,Q) = \frac{1}{L} \sum_{\ell=1}^L W_p^p(P_{\theta_\ell}, Q_{\theta_\ell}),
\]
with \(O(L^{-1/2})\) convergence~\cite{nadjahi2020statistical}, improvable to \(O(L^{-3/4})\) with spherical harmonics control variates (SHCV)~\cite{nguyen2023control}. However, unlike block SLQ variants, SW methods degrade in high dimensions and do not exploit matrix sparsity or structure.

\subsection{Notation} \label{sec:notation}
Vectors are lowercase (e.g., \(x, y\)); matrices are uppercase (e.g., \(A, \Sigma\)), with transpose \(A^T\). For \(A \in \mathbb{S}_+^n\), let \(\lambda_i(A)\) be its eigenvalues in decreasing order. We write \(\|A\|\), \(\|A\|_F\), and \(\tr(A)\) for the spectral norm, Frobenius norm, and trace, respectively. \(\mathbb{S}_+^n\) denotes the cone of symmetric positive semidefinite \(n\times n\) matrices. \(\mathcal{N}(\mu, \Sigma)\) is the Gaussian with mean \(\mu\) and covariance \(\Sigma\). For \(S \subset \{1,\dots,n\}\), \(A_S = A(S,S)\) is the submatrix of \(A\) indexed by \(S\), and \(P_S\) projects onto \(\{e_i: i\in S\}\). We write \(A \preceq B\) to mean \(B - A\) is positive semidefinite.

Expectations and variances are denoted \(\mathbb{E}[\cdot]\), \(\operatorname{Var}[\cdot]\), and \(\|x\|\) is the Euclidean norm.


\subsection{KL Divergence of Gaussians as a Trace Problem}

\noindent
The KL divergence between Gaussians \(p = \mathcal{N}(\mu_1, \Sigma_1)\) and \(q = \mathcal{N}(\mu_2, \Sigma_2)\), with \(\mu_1 = \mu_2 \in \mathbb{R}^d\) and \(\Sigma_1, \Sigma_2 \in \mathbb{S}_+^d\), simplifies to \(D_{\mathrm{KL}}(p\parallel q) = \frac{1}{2}[\tr(\Sigma_2^{-1}\Sigma_1) + \tr(\log\Sigma_1) - \tr(\log\Sigma_2) - d]\)~\cite{abdar2021review, gultekin2017nonlinear, furrer2007estimation, ledoit2004well, fan2008covariance, srivastava2005tests}. Letting \(\Sigma_2^{-1} = LL^T\) (e.g., from PDE-generated covariances or Vecchia approximations~\cite{fuentes2023estimating,furrer2007estimation}) and \(\Sigma = \Sigma_1\), this becomes \(D_{\mathrm{KL}} = \frac{1}{2} \tr(f(LL^T\Sigma))\) where \(f(\lambda) = \lambda - \log\lambda - 1\). To approximate \(\tr(f(A))\), we use Hutchinson’s method~\cite{hutchinson1989stochastic} with i.i.d.\ Gaussian probes \(g_i \sim \mathcal{N}(0,I)\): 
\[\tr(f(A)) \approx \frac{1}{q} \sum_{i=1}^q g_i^T f(A) g_i.\] 
Applying Lanczos quadrature~\cite{ubaru2017fast, bai1996some} to each \(g_i\), we construct a tridiagonal matrix \(T_i\), diagonalize \(T_i = U_i \Lambda_i U_i^T\), and obtain \(\sum_{j=1}^k w_{ij} f(\mu_{ij})\) where \(w_{ij} = (U_i)_{1j}^2\). Our estimator is
\[
\widehat{D}_{\mathrm{KL}}(f(A)) := \frac{1}{2q} \sum_{i=1}^q \sum_{j=1}^k w_{ij} f(\mu_{ij}).
\]




\subsection{Wasserstein-2 Distance as a Trace Problem}

\noindent
For Gaussians \(p = \mathcal{N}(\mu, \Sigma_1)\), \(q = \mathcal{N}(\mu, \Sigma_2)\), the squared Wasserstein-2 distance is given by $W_2^2(p, q) = \tr(\Sigma_1 + \Sigma_2) - 2\, \tau$ where $\tau:=\tr((\Sigma_1^{1/2} \Sigma_2 \Sigma_1^{1/2})^{1/2})$. Noting \(\spec(\Sigma_1^{1/2} \Sigma_2 \Sigma_1^{1/2}) = \spec(\Sigma_2 \Sigma_1)\) under the assumption that \(\Sigma_1 \in \mathbb{S}_+^n\), since \(\Sigma_2 \Sigma_1 = \Sigma_1^{-1/2} (\Sigma_1^{1/2} \Sigma_2 \Sigma_1^{1/2}) \Sigma_1^{1/2}\), define \(B = \Sigma_2 \Sigma_1\), so that \(\tr((\Sigma_1^{1/2} \Sigma_2 \Sigma_1^{1/2})^{1/2}) = \tr(f(B))\) with \(f(x) = \sqrt{x}\).

Applying stochastic Lanczos quadrature on \(B\) (no explicit square roots required): for \(i = 1, \dots, q\), draw \(g_i\), run \(k\) Lanczos steps to get the tridiagonal \(T_i = U_i \Theta_i U_i^T\), where \(\Theta_i = \mathrm{diag}(\theta_{i1}, \dots, \theta_{ik})\), and compute \(g_i^T f(B) g_i \approx e_1^T f(T_i) e_1 = \sum_{j=1}^k w_{ij} \sqrt{\theta_{ij}}\), with \(w_{ij} = (U_i)_{1j}^2\). Then \(\widehat{\tau} = \frac{1}{q} \sum_{i=1}^q \sum_{j=1}^k w_{ij} \sqrt{\theta_{ij}}\), and \(W_2^2(p, q) \approx \tr(\Sigma_1 + \Sigma_2) - 2 \widehat{\tau}\). The trace term \(\tr(\Sigma_1 + \Sigma_2)\) can also be efficiently estimated using BOLT, requiring only matrix–vector products with \(\Sigma_1\) and \(\Sigma_2\).


\section{Block-Orthonormal Stochastic Lanczos Quadrature} \label{sec:blockslq}

\noindent
We introduce the Block-Orthonormal Stochastic Lanczos Quadrature (BOLT) method, a randomized trace estimator that avoids the high overhead of XTrace~\cite{epperly2024xtrace} and the projection step used in Hutch++~\cite{meyer2021hutch++}. While Hutch++ requires a randomized SVD and assumes fast spectral decay, BOLT applies Lanczos directly to orthonormal blocks, achieving the same convergence rate without low-rank assumptions. A related idea of using designed probing for Gaussian processes was explored in \cite{stein2013stochastic}, which observed improved empirical performance but did not analyze the asymptotic behavior of the estimator.  We present BOLT in Algorithm~\ref{algo2} and compare its FLOP costs to other stochastic trace estimators in detail in Appendix~\ref{appendix:trace_flops}. All code to reproduce the figures can be found at \url{https://github.com/chebyshevtech/BOLT}.

\begin{algorithm}[h]
\caption{Block-Orthonormal SLQ (BOLT): Randomized Estimation of \(\tr(f(A))\)}
\label{algo2}
\begin{algorithmic}[1]
\STATE \textbf{Input:} Symmetric matrix or operator \(A\), analytic function \(f\), number of probes \(q\), Lanczos steps \(k\), block size \(b\).
\STATE \textbf{(KL case)}: set \(A = L^T\Sigma L\), \(f(\lambda)=\lambda - \ln\lambda - 1\).
\FOR{\(i = 1,\dots,q\)}
  \STATE Draw \(Z_i\in\mathbb{R}^{n\times b}\) with either Rademacher or standard Gaussian entries, and orthonormalize to obtain \(V_i\).
  \STATE Run block Lanczos on \(A\) starting from \(V_i\) for \(k\) steps to form tridiagonal \(T_i\).
  \STATE Diagonalize \(T_i = U_i\,\mathrm{diag}(\mu_{i1},\dots,\mu_{ikb})\,U_i^T\).
  \STATE Compute weights \(w_{ij} = \sum_{r=1}^b (U_i)_{rj}^2\) and 
  \(\eta_i = \sum_{j=1}^{kb} w_{ij}\,f(\mu_{ij})\).
\ENDFOR
\STATE \(\displaystyle\widehat{\tr}(f(A)) = \frac{n}{q b}\sum_{i=1}^q \eta_i.\)
\STATE \textbf{Output:} \(\widehat{\tr}(f(A))\) (or for KL, \(\widehat D_{\mathrm{KL}}=\tfrac12\,\widehat{\tr}(f(A))\)).
\end{algorithmic}
\end{algorithm}

The SLQ method which was \cite{bai1996some} combines Hutchinson probing with Gaussian quadrature. Rather than computing the full matrix at \(O(n^3)\) cost where $n$ denotes the size of the matrix, SLQ uses \(qk\) matrix–vector products \cite{kuczynski1992estimating} for certain choice of $q$ and $k$. It was later shown that for any fixed \(q\),  \(k = O(\ln n / \sqrt{\epsilon})\) Lanczos steps can yield a precision of order \(\epsilon\). Our approach, motivated by block SLQ based method of~\cite{fuentes2023estimating} introduces a block-orthonormal SLQ variant for improved convergence. 

To rigorously quantify the performance of our estimators, we now present formal theoretical results establishing unbiasedness, variance bounds, and exact recovery guarantees under appropriate sampling conditions. These results not only underpin the empirical observations reported in later sections but also highlight the fundamental differences between block and scalar probing strategies. Full proofs of the following theorems are deferred to Appendix~\ref{appendix:proofs}, specifically Sections~\ref{appendix:blockSLQ}.

\begin{theorem}[Block SLQ Estimator] \label{thm:blockSLQ}
Let \(A \in \mathbb{R}^{n\times n}\) be symmetric positive semidefinite with eigenvalues \(\{\lambda_i\}_{i=1}^n\), and let
\(f:[0,\infty)\to\mathbb{R}\)
be convex and twice continuously differentiable. For each of \(q\) independent trials, generate an \(n\times b\) random matrix \(V\) with orthonormal columns (\(\mathbb{E}[VV^T]=\frac{b}{n}I_n\)) and define
\(
X(V) := \frac{n}{b}\,\operatorname{tr}\bigl(V^T f(A) V\bigr).
\)
The Monte Carlo estimator
\[
\widehat{T} := \frac{1}{q} \sum_{i=1}^q X(V_i)
\]
satisfies:
\begin{enumerate}
\item \(\mathbb{E}[X(V)] = \operatorname{tr}\bigl(f(A)\bigr)\). (If \(b=n\), then \(VV^T = I_n\) and \(X(V)\) is exact.)

\item The variance is
\[
\operatorname{Var}[X(V)] = \frac{2n}{b(n+2)}\left(1-\frac{b-1}{n-1}\right) \left( \sum_{i=1}^n f(\lambda_i)^2 - \frac{1}{n}\left(\operatorname{tr}(f(A))\right)^2 \right),
\]
and
\[
\operatorname{Var}[\widehat{T}] = \frac{1}{q} \operatorname{Var}[X(V)].
\]

\item If \(X(V)\in[\ell_{\min},\ell_{\max}]\) almost surely, then for any \(\varepsilon>0\),
\[
\mathbb{P}\left(\bigl|\widehat{T}-\operatorname{tr}(f(A))\bigr| \geq \varepsilon\right) \leq 2\exp\left(-\frac{2q\varepsilon^2}{(\ell_{\max}-\ell_{\min})^2}\right),
\]
so that with probability at least \(1-\delta\),
\[
\bigl|\widehat{T}-\operatorname{tr}(f(A))\bigr| \leq (\ell_{\max}-\ell_{\min})\sqrt{\frac{\ln(2/\delta)}{2q}}.
\]
\end{enumerate}
\end{theorem}

\medskip

\begin{figure}[h]
    \centering
    \includegraphics[width=0.5\textwidth]{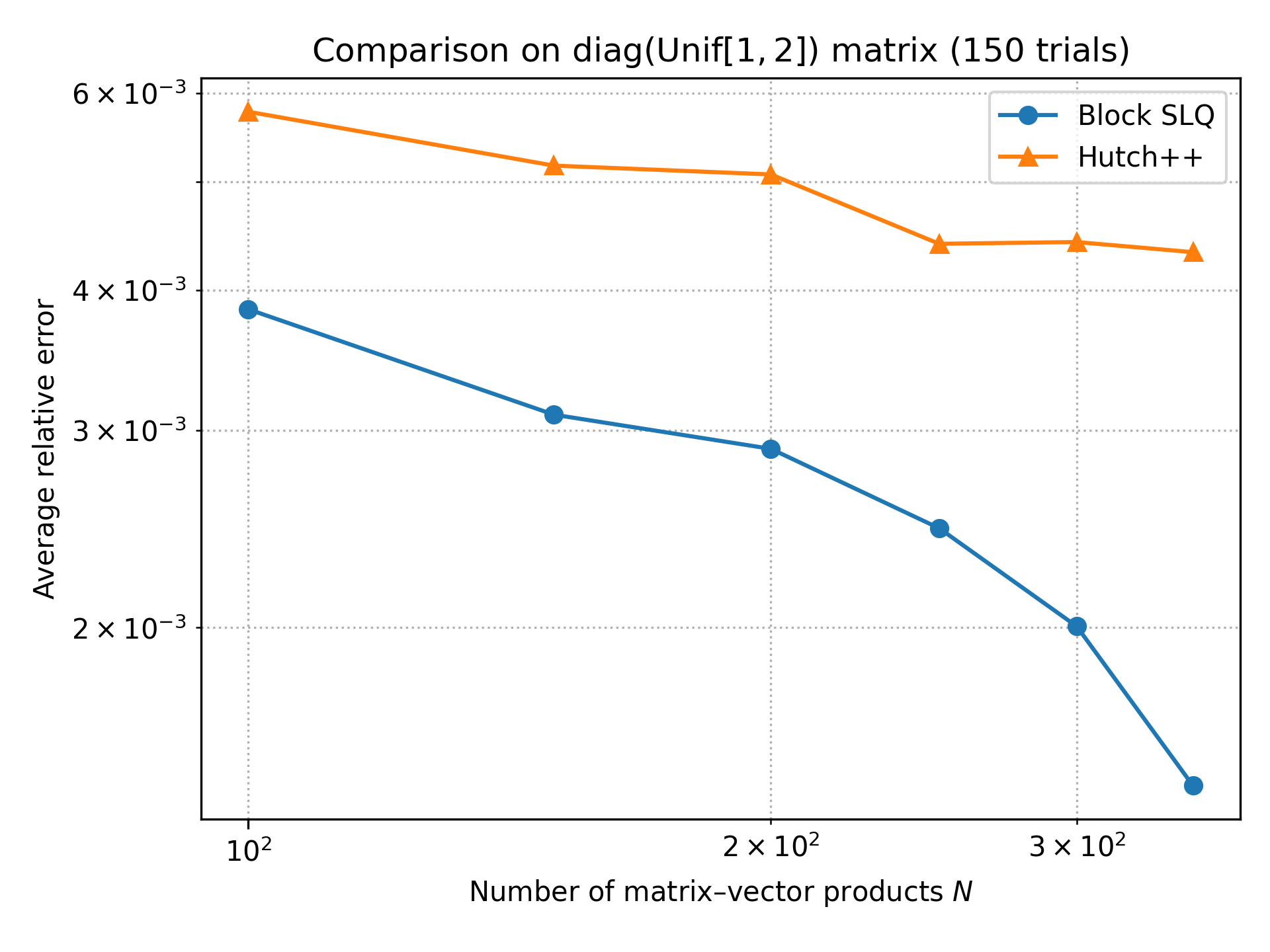}
    \caption{BOLT outperforms Hutch++ in approximating the trace of \(\mathrm{diag}(\mathrm{Unif}[1,2])^2\) (150 trials)}
    \label{fig:hutchppfail}
\end{figure}

A key advantage of block–orthonormal SLQ arises directly from its variance.  
As shown in the proof of Theorem~\ref{thm:blockSLQ}, BOLT achieves strictly lower
variance than Hutchinson for every block size \(b>1\).  
The comparison with Hutch++ is governed by the variance ratio derived in the same
proof and stated explicitly in \eqref{eq:bolt-vs-hpp-var}, which we reproduce here for
convenience:
\begin{equation}
\frac{\operatorname{Var}[X(V)]}{\operatorname{Var}[\widehat T_{\mathrm{H++}}(A)]}
\le
\frac{mn}{3b(n+2)}
\left(1 - \frac{b-1}{n-1} \right)
\cdot
\frac{\sum_{i=1}^n \lambda_i^2 - \frac{1}{n} \left( \sum_{i=1}^n \lambda_i \right)^2}
     {\sum_{i=k+1}^n \lambda_i^2}.
\end{equation}
BOLT can outperform Hutch++ in regimes where the spectrum is nearly flat and the
residual tail energy \(\sum_{i=k+1}^n \lambda_i^2\) remains large.  
For such flat-spectrum matrices—meaning that for all eigenvalues
\(\lambda_i, \lambda_j \in \mathrm{Spec}(A)\),
\[
|\lambda_i - \lambda_j| < d, \qquad 0 \le d \le \varepsilon n,
\]
the variance ratio becomes particularly favorable, leading to rapid concentration,
consistent with the empirical behavior in Figure~\ref{fig:hutchppfail}.

\begin{remark}
It is worth noting that if one uses the block size \(b=1\), then Block SLQ reduces to
the scalar SLQ Hutchinson estimator of \cite{bai1996some}, which is known to converge
at rate \(O(N_{\rm mv}^{-1/2})\). Likewise, if one chooses the maximal block size
\(b=n\), then the estimator becomes deterministic, since the orthonormal probe spans
the entire space and \(V V^{\top} = I_n\), yielding the exact trace.

\end{remark}


\begin{figure}[htbp]
    \centering
    \begin{minipage}[t]{0.48\textwidth}
        \centering
        \includegraphics[width=\textwidth]{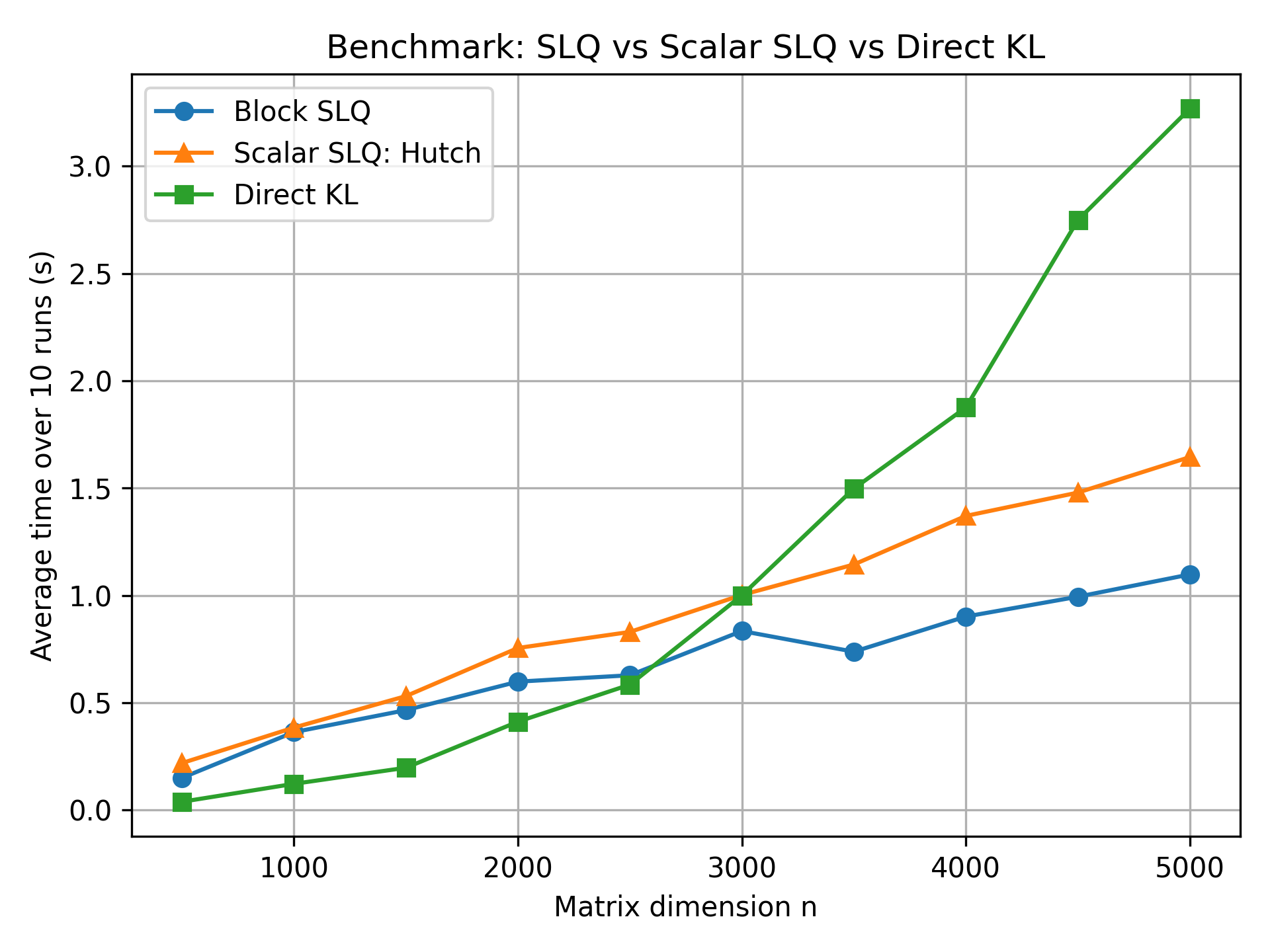}
        \caption{Timing comparisons}
        \label{fig:timing_comparison}
    \end{minipage}
    \hfill
    \begin{minipage}[t]{0.48\textwidth}
        \centering
        \includegraphics[width=\textwidth]{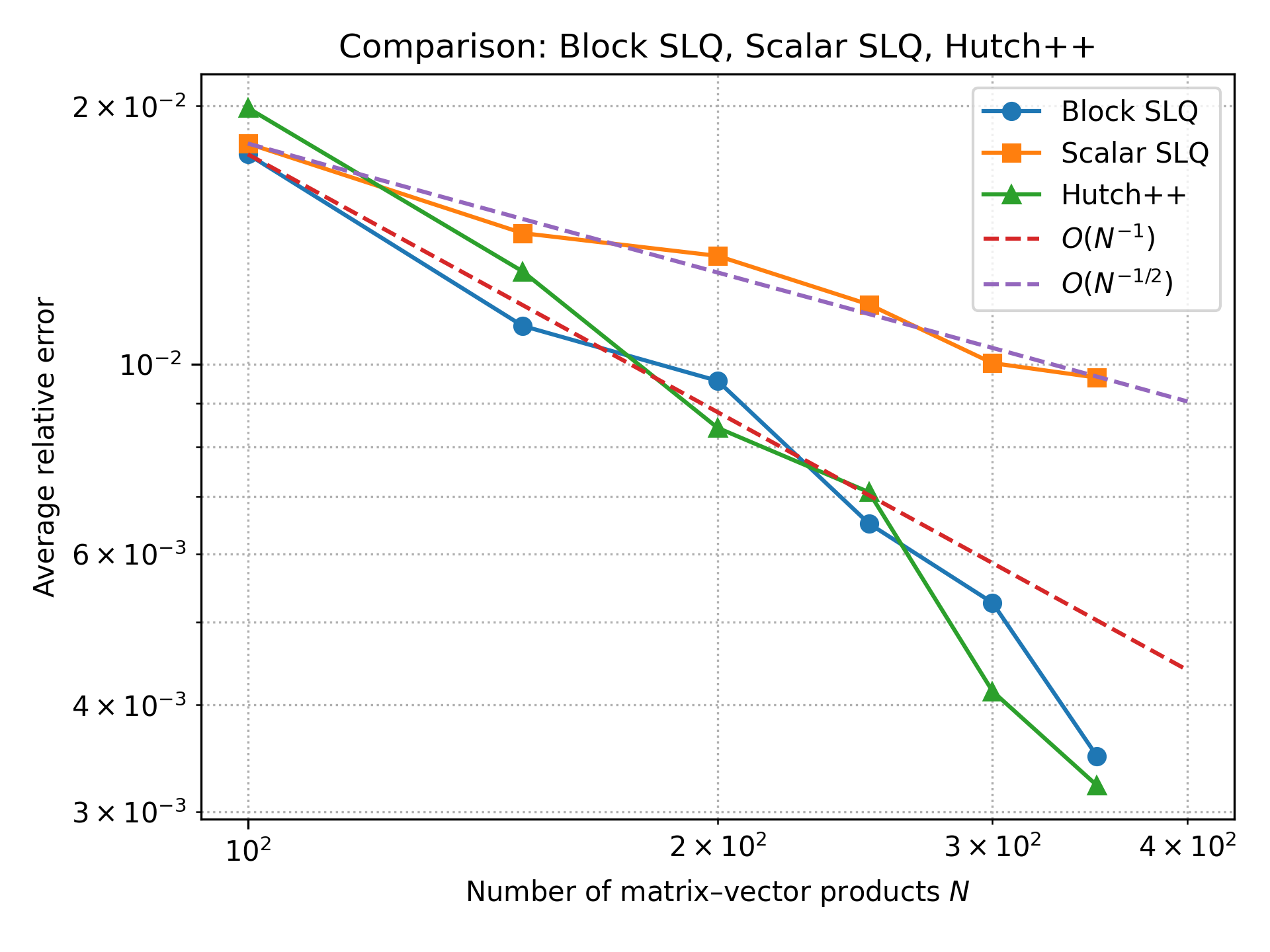}
        \caption{Relative error (KL divergence) vs.\ number of matrix–vector products (Scalar-Hutch vs.\ Block SLQ) vs. Hutch++}
        \label{fig:error_vs_matvec_actual_N}
    \end{minipage}
\end{figure}



We now numerically demonstrate KL divergence estimation between two zero-mean
Gaussians, where \(A = \Sigma \in \mathbb{R}^{200\times200}\) is a Gaussian RBF kernel with
\(\sigma = 2\), and \(LL^T\) is the Cholesky factor of the precision matrix.
Figure~\ref{fig:timing_comparison} compares wall-clock runtimes: SLQ- and
Hutchinson-based estimators scale efficiently with matrix size, unlike direct log-det
computation. Figure~\ref{fig:error_vs_matvec_actual_N} illustrates the convergence
behavior of the estimators: Block SLQ consistently outperforms Hutchinson due to its
strictly lower variance, and matches
that of Hutch++, reflecting the variance advantage discussed earlier.

Unlike Hutch++ \cite{meyer2021hutch++}, which first multiplies the matrix of interest $A$ by a sketching matrix $S$ (a randomized matrix with independent Gaussian or Rademacher entries) to form a sketch $AS$ and then performs a QR-decompsition of $AS$ to approximate the top eigenspace, our block SLQ method takes a simpler route: we draw a single Rademacher or Gaussian block $Z \in \mathbb{R}^{n \times b}$, compute its QR factorization $Z = \widetilde Q\,\widetilde R$, and use the orthonormal columns of $\widetilde Q$ directly as our block probe—without any prior multiplication by $A$. This minor change still induces negative covariances between the quadratic forms $v_r^T f(A)\,v_r$, since orthogonal directions naturally “self-correct”: overestimation by one vector implies underestimation by its orthogonal complement.

Figure~\ref{fig:hutchppfail} shows this effect for \(A = \mathrm{diag}(d)\) with \(d_i \sim \mathrm{Unif}[1,2]\). Hutch++ builds a sketch \(Q \in \mathbb{R}^{n \times s}\) with \(s = N_{\mathrm{mv}}/3\), computes \(\tr(Q^T A Q)\) exactly, and estimates the remainder with \(g = N_{\mathrm{mv}} - s\) Monte Carlo probes, yielding \(O(g^{-1/2}) = O(N_{\mathrm{mv}}^{-1/2})\) error. Block SLQ instead draws one orthonormal block \(\widetilde Q \in \mathbb{R}^{n \times b}\) and computes $\widehat T = \frac{n}{b} \sum_{j=1}^b d_{k_j}.$

\section{Sub-block Stochastic Lanczos Quadrature} \label{sec:subblockslq}

The goal of this section is address the trace estimation of matrix functions when the matrix vector product computation of the full matrix is unavailable or too expensive. We seek to propose a consistent trace estimation method which works under the assumption when the matrix vector product sub-blocks are available. Using the the following result (proved in Appendix~\ref{appendix:unbiasedsubblockbound}), we propose a unified estimator (Algorithm~\ref{algo:unified}) that applies even when \(A = L^T\Sigma L\) is singular or \(\Sigma\) is large and sparse, relying on small submatrices instead of full matrix–vector products to approximate \(\tr(f(A))\).

We now formalize the correctness of the sub-block based trace estimation.

\begin{lemma}[Unbiased Subblock Trace Estimation] \label{lem:unbiasedsubblockbound}
Let \(A\in\mathbb{R}^{n\times n}\) be symmetric positive semidefinite and \(f\colon[0,\infty)\to\mathbb{R}\) nondecreasing. For any subset \(S\subset\{1,\dots,n\}\) with \(|S|=s\), define
\[
X(S) := \frac{n}{s}\,\tr\bigl(f(A)_S\bigr).
\]
See section \ref{sec:notation} for $f(A)_S$. Then \(\mathbb{E}_S[X(S)] = \tr(f(A))\)
\end{lemma}

\begin{figure}[h]
    \centering
    \begin{minipage}{0.48\textwidth}
        \centering
        \includegraphics[width=\textwidth]{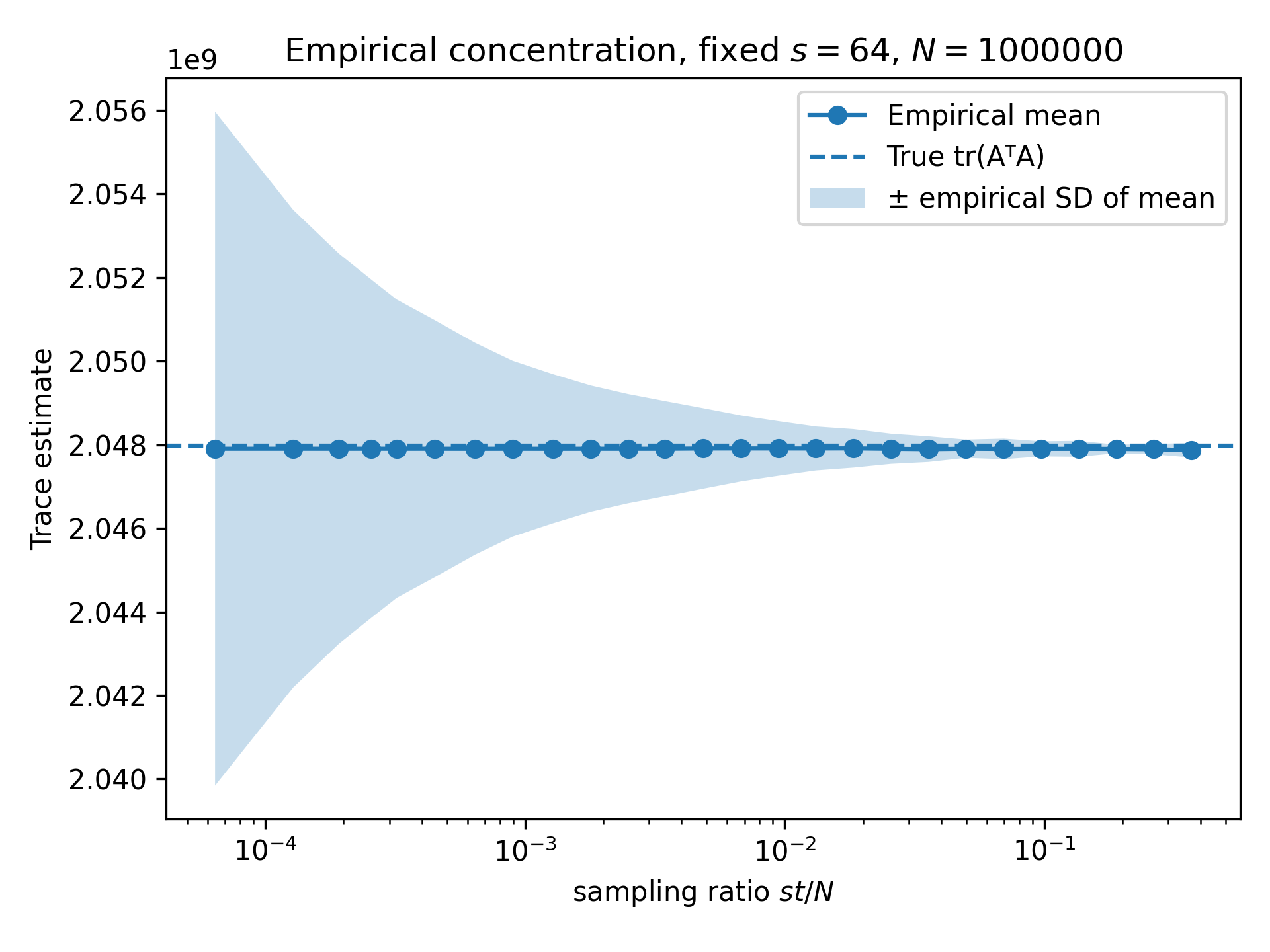}
        \caption{Empirical trace recovery versus sampling ratio \(st/N\).}
        \label{fig:trace_recovery}
    \end{minipage}
    \hfill
    \begin{minipage}{0.48\textwidth}
        \centering
        \includegraphics[width=\textwidth]{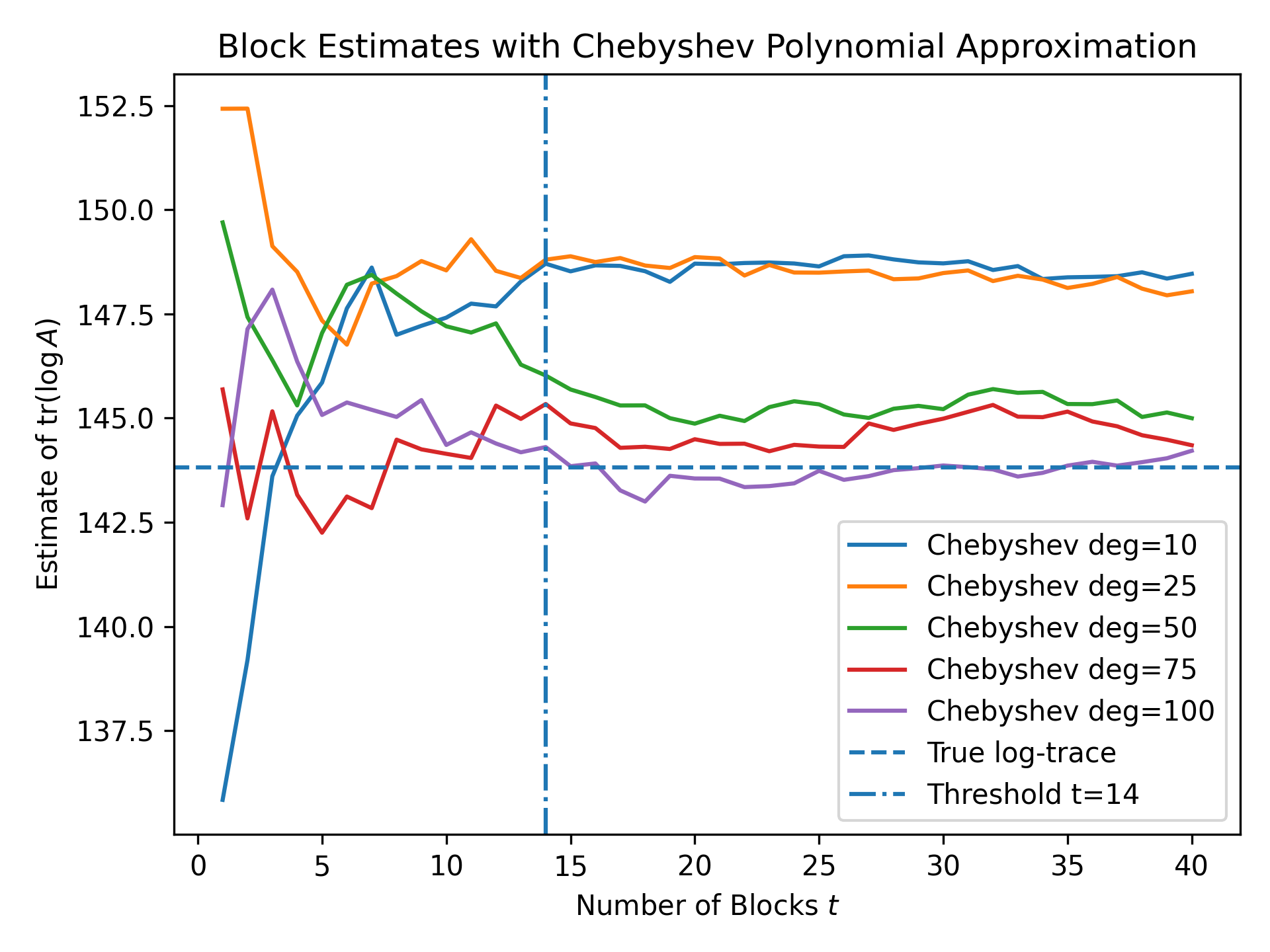}
        \caption{Polynomial localization via Chebyshev block estimates.}
        \label{fig:chebyshev_block_estimates}
    \end{minipage}
\end{figure}

Figure~\ref{fig:trace_recovery} reports an empirical study of trace recovery using the unbiased subblock estimator from Lemma~\ref{lem:unbiasedsubblockbound}. 
We generate a random matrix \(B \in \mathbb{R}^{m \times n}\) with independent standard normal entries, with \(m=2048\) and \(n=10^6\), and form \(A = B^T B\), so that \(\tr(A) = \|B\|_F^2\).
For each block, we sample a subset \(S \subset \{1,\dots,n\}\) of size \(s=64\) uniformly and compute the scaled subblock trace
\[
X(S) = \frac{n}{s} \sum_{i \in S} A_{ii},
\]
which is an unbiased estimator of \(\tr(A)\).
We generate many independent block estimates and study the concentration of the averaged estimator
\(
\frac{1}{t}\sum_{k=1}^t X(S_k)
\)
as the number of blocks \(t\) increases.
Uncertainty is quantified empirically by forming multiple independent \(t\)-block estimators and reporting the resulting standard deviation.
The horizontal axis shows the normalized sampling ratio \(st/n\), showing that accurate trace recovery is achieved while observing only a vanishing fraction of the diagonal entries.

When \(A\) is too large or singular to apply \(f(A)\) directly, we instead approximate \(f(A_S)\) using local computations. The next result (proved in~\ref{appendix:localization}) ensures this approximation is valid for polynomial filters:

\begin{theorem}[Exact Localization of Polynomial Filters on Principal Subblocks] \label{thm:localization}
Let \(A\in\mathbb{R}^{n\times n}\) have sparsity pattern inducing a graph \(G\), and assume \(\mathrm{spec}(A)\subset[\lambda_{\min},\lambda_{\max}]\). Let \(p_m\) be a degree-\(m\) Chebyshev polynomial approximating \(f\) with uniform error \(\varepsilon\). For any subset \(S\), define
\[
S_r := \{\,v : \mathrm{dist}_G(v,S) \le r\,\}.
\]
If \(r\ge m\), then
\[
[p_m(A)]_{S,S} = [p_m(A_{S_r,S_r})]_{S,S} = p_m(A_S).
\]
Moreover, the trace approximation error on \(S\) satisfies
\[
\bigl|\tr\bigl(f(A)_S\bigr) - \tr\bigl(f(A_S)\bigr)\bigr|\;\le\;2s\,\varepsilon.
\]
\end{theorem}

\begin{remark}
In practice, we do not construct \(p_m\) explicitly. Instead, block Lanczos or Gaussian quadrature implicitly builds such filters and integrates polynomials exactly up to degree \(2k-1\) with \(k\) iterations.
\end{remark}

Figure~\ref{fig:chebyshev_block_estimates} confirms this localization. Such localization techniques have been widely employed in convolutional neural networks for graph filtering~\cite{defferrard2016convolutional}. An alternative approach to stochastic Lanczos quadrature is to approximate the matrix logarithm using a deep polynomial representation~\cite{wei2025deep}, designed to minimize the number of matrix multiplications needed for accurate reconstruction~\cite{jarlebring2024polynomial}.

\begin{algorithm}[h]
\caption{Unified Block SLQ: Randomized Estimation of \(\tr(f(A))\)}
\label{algo:unified}
\begin{algorithmic}[1]
\STATE \textbf{Input:} Symmetric matrix/operator \(A\), analytic \(f\), probes \(q\), subblocks \(t\), block size \(s\), steps \(k=s\), tol.\ \(\varepsilon\).
\STATE Compute effective dimension \(r_{\mathrm{eff}}=\#\{i:A_{ii}>\varepsilon\}\).
\IF{\(t=1\) or \(r_{\mathrm{eff}}\le s\)} \STATE Reduce to Algorithm~\ref{algo2} \ENDIF
\STATE Sample subblocks \(S_1,\dots,S_t\subset\{1,\dots,r_{\mathrm{eff}}\}\), each of size \(s\).
\FOR{\(i=1,\dots,t\)}
  \STATE Form \(A_{S_i}=A(S_i,S_i)\).
  \FOR{\(j=1,\dots,q\)}
    \STATE Draw \(Z_{ij}\in\mathbb{R}^{s\times b}\) with either Rademacher entries (\(\pm1\)) or standard Gaussian entries, and orthonormalize to obtain \(V_{ij}\).
    \STATE Run block Lanczos on \(A_{S_i}\) from \(V_{ij}\) for \(k\) steps to get \(T_{ij}\).
    \STATE Diagonalize \(T_{ij}=U_{ij}\,\mathrm{diag}(\mu_{ij1},\dots,\mu_{ij,ks})\,U_{ij}^T\).
    \STATE \(\eta_{ij}=\sum_{l=1}^{ks}\bigl(\sum_{r=1}^s (U_{ij})_{rl}^2\bigr)\,f(\mu_{ijl})\).
  \ENDFOR
  \STATE \(\eta_i=\frac1q\sum_{j=1}^q\eta_{ij}.\)
\ENDFOR
\STATE \(\widehat{\tr}(f(A))=\frac{r_{\mathrm{eff}}}{ts}\sum_{i=1}^t\eta_i.\)
\STATE \textbf{Output:} \(\widehat{\tr}(f(A))\) (for KL, \(\widehat D_{\mathrm{KL}}=\tfrac12\widehat{\tr}(f(A))\)).
\end{algorithmic}
\end{algorithm}

\subsection{Applications of Proxy KL Divergence}

In this subsection, we will define the proxy KL divergence and explore its applicability in two ways. First, we verify that in the low sample regime, sparsity in the Cholesky factor of the precision matrix is necessary to produce accurate results. Second, we use our proxy KL as a regularizer for the classification of MNIST data. In both cases, the matrix of interest is singular; therefore, the regular KL divergence is not defined.
One such instance arises for the sample covariance \(\widetilde{\Sigma} = \sum_{i=1}^m u_i u_i^T\), with \(u_i \sim \mathcal{N}(0,\Sigma)\) and $m\ll n$. Notice that $\widetilde{\Sigma}$ is singular, even though the true covariance \(\Sigma\) is full rank. This makes \(\log(\widetilde{\Sigma})\) and the standard KL divergence ill-defined.
We introduce a \emph{proxy KL divergence} via trace functionals on subblocks of \(\widetilde{\Sigma}\), which remains well-posed in singular, sample-starved settings and applies to any symmetric matrix, regardless of rank or access.

\begin{lemma}[Subblock full‐rankness from Wishart sampling] \label{thm:fullrank}
Let \(\Sigma\in\mathbb R^{n\times n}\) be positive definite, and let $u_1, \dots, u_m \sim \text{iid}~\mathcal{N}(0, \Sigma)$
  Define the sample covariance (unnormalized)
$
\widetilde\Sigma \;=\;\sum_{i=1}^m u_i\,u_i^T,
$
which has rank at most \(m\).  Fix any index set \(S\subset\{1,\dots,n\}\) of size \(s\le m\), and let \(P_S\in\{0,1\}^{n\times s}\) be the coordinate‐selection matrix picking those \(s\) rows.  Then the \(s\times s\) subblock
$\widetilde\Sigma_S \;=\;P_S^T\,\widetilde\Sigma\,P_S$
follows the central Wishart distribution \(W_s(P_S^T\Sigma P_S,\;m)\).  Since \(m\ge s\) and \(P_S^T\Sigma P_S\) is invertible, \(\widetilde\Sigma_S\) is almost surely invertible and hence \(\text{rank}(\widetilde\Sigma_S)=s\).
\end{lemma}

Although \(\widetilde{\Sigma} = \sum_{i=1}^m u_i u_i^T\) is rank-deficient when \(m < n\), its subblocks \(\widetilde{\Sigma}_S\) are full rank with probability one whenever \(s \le m\), a property inherited from the Wishart structure (Lemma~\ref{thm:fullrank}, proved in~\ref{appendix:fullrank}). For stable KL estimation, one must choose subblocks with \(s \le m\); beyond this threshold, \(\widetilde{\Sigma}_S\) becomes singular (Figure~\ref{fig:hist2}). Appendix~\ref{appendix:eigenvalue_distributions} discusses this transition further, including smallest eigenvalue statistics from Edelman~\cite{edelman1988eigenvalues} and universality extensions by Tao and Vu~\cite{tao2010random}.

\begin{definition}[Proxy KL Divergence]
Let \(A\in\mathbb{R}^{n\times n}\) be symmetric positive semidefinite (possibly singular) and define 
\[
f(\lambda) := \lambda - \log\lambda - 1.
\]
For any subset \(S\subset\{1,\dots,n\}\) with \(|S|=s\), let \(f(A)_S\) denote the \(s\times s\) subblock of \(f(A)\) indexed by \(S\).  The \emph{proxy KL divergence} is then
\[
\widehat{D}_{\mathrm{KL}}^{\mathrm{proxy}}(f(A))
\;:=\;  \frac{1}{t}\sum_{i=1}^t \widehat{D}_{\mathrm{KL}} (f(A)_{S_i}) 
\]
where $S_1, \ldots , S_t$ are independent sub-blocks sampled randomly from the set $\{1,\ldots , n\}$ without replacement. 

\end{definition}

Notice that \(\mathbb{E}\bigl[\tfrac{n}{s}\tr(f(A)_S)\bigr]=\tr(f(A))\) by Lemma~\ref{lem:unbiasedsubblockbound}. Hence when the matrix is full rank, \(\widehat{D}_{\mathrm{KL}}^{\mathrm{proxy}}\) is an unbiased estimator of the true KL divergence \(\tfrac12\tr(f(A))\).

We illustrate our algorithm with two examples. First, we apply Algorithm~\ref{algo:unified} to estimate both the standard KL divergence in the full-rank case and a proxy KL divergence when the matrix is singular. Figures~\ref{fig:err} and~\ref{fig:err_sparse} in appendix~\ref{appendix:figures} show this behavior under least-squares and sparsity-constrained Cholesky factorizations. Without sparsity in \(L\), the error matrix \(E = LL^T\Sigma\) has \(\|E\|_F = 2\times10^{14}\), matching the proxy KL divergence and indicating poor approximation. Enforcing sparsity in \(L\) sharply reduces this error, as reflected in the KL metric.

\begin{figure}[h]
    \centering
    \includegraphics[width=1\textwidth]{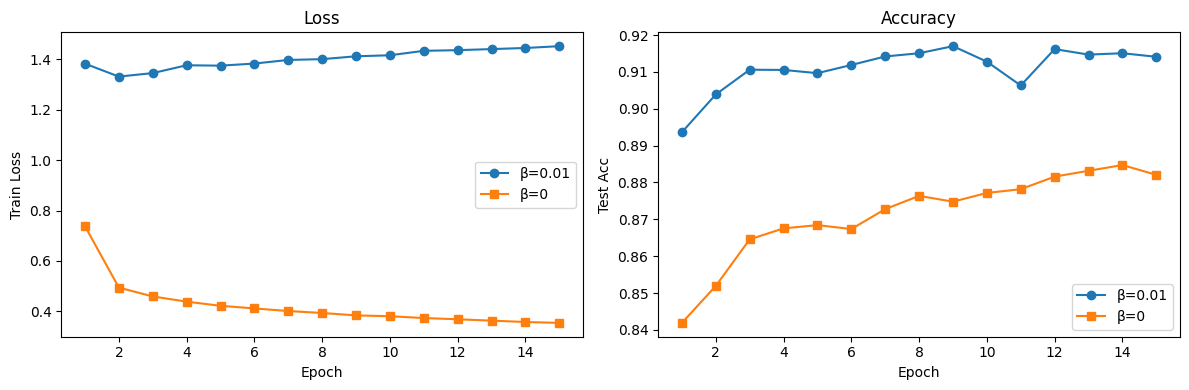}
    \caption{Test accuracy per epoch comparing models with and without SLQ-KL regularization. Regularization (\(\beta=0.01\)) improves generalization accuracy.}
    \label{fig:accuracy_sql}
\end{figure}

Second, we train an MLP (input 784 $\rightarrow$ hidden 8 $\rightarrow$ output 10) on MNIST \cite{deng2012mnist} with batch size 4. At each step, we compute hidden activations \(H \in \mathbb{R}^{4\times 8}\), center them to \(H_c\), and form the empirical covariance \(\Sigma = \frac{1}{3} H_c^T H_c \in \mathbb{R}^{8\times 8}\). Since \(\Sigma\) is singular, standard KL regularization is undefined. We instead apply the SLQ-KL penalty \(\frac{\beta}{2} \tr(f(\Sigma))\) via subblock SLQ to encourage isotropy, a property known to aid generalization in NLP and vision~\cite{arora2017simple, mu2017all, cai2021isotropy}. Figure~\ref{fig:accuracy_sql} shows improved test accuracy with \(\beta=0.01\); Figure~\ref{fig:labels_pred} highlights corrected misclassifications. See Appendix~\ref{appendix:experiment-details} for full setup.

\begin{figure}[h]
    \centering
    \includegraphics[width=1\textwidth]{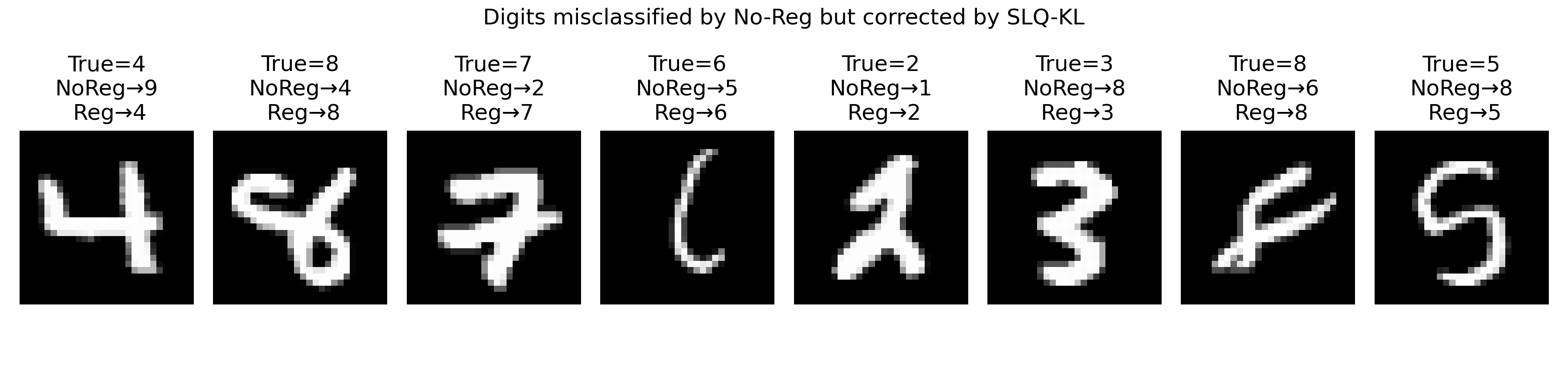}
    \caption{MNIST test examples where the SLQ-KL regularized model (Reg) correctly classifies digits that the unregularized model (NoReg) misclassifies. Top: true labels. Bottom: predicted labels for NoReg and Reg.}
    \label{fig:labels_pred}
\end{figure}

\subsection{Application to Divergence Estimation for Hierarchical Off-Diagonal Low-Rank (HODLR) Matrices}
\label{sec:hodlr}

In this subsection, we introduce an efficient method for computing statistical divergences of hierarchical off-diagonal low-rank (HODLR) matrices using only matrix–vector products.

Hierarchical off-diagonal low-rank (HODLR) matrices provide an efficient multilevel
representation for dense yet structured operators such as covariance and kernel
matrices. Although HODLR factorizations can be constructed for more general
matrices, here we introduce a symmetric variant tailored to covariance settings. A
symmetric matrix \(A \in \mathbb{R}^{n\times n}\) can be recursively partitioned into
\(2\times 2\) blocks,
\[
A = 
\begin{bmatrix}
A_{11} & A_{12} \\
A_{21} & A_{22}
\end{bmatrix},
\qquad
A_{12},\,A_{21}^\top \text{ are low rank.}
\]
At level \(\ell\), each off-diagonal block \(A_{ij}^{(\ell)}\) is approximated by a factorization
\[
A_{ij}^{(\ell)} \approx U_{ij}^{(\ell)} S_{ij}^{(\ell)} {V_{ij}^{(\ell)}}^{\!\top},
\qquad
\operatorname{rank}(S_{ij}^{(\ell)}) = r_\ell \ll n_\ell.
\]
This yields a hierarchical factorization of \(A\) requiring only \(O(n\log n)\) storage and allowing matrix–vector products to be performed in \(O(n\log n)\) operations.

Following Lin et al. \emph{peeling algorithm} \cite{lin2011fast}, the HODLR factors can be constructed directly from matrix–vector products, without ever forming the full matrix. 
Let \(A:\mathbb{R}^n \to \mathbb{R}^n\) denote the operator accessible through mat–vec evaluations.
At each level, we recursively subtract the already recovered low-rank interactions from the residual operator and apply randomized SVD to the remaining off-diagonal blocks.
With \(k\) levels and target ranks \(\{r_\ell\}\), the total complexity of building the HODLR representation is
\[
O\!\left( \sum_{\ell=1}^k r_\ell n_\ell \right) = O(n\log n),
\]
assuming uniform block sizes and bounded ranks.

Although the peeling algorithm provides a provably accurate reconstruction of the off-diagonal blocks, there is currently no fast way to evaluate the global approximation error.  
For example, computing the Frobenius norm
\(\|A - \tilde A\|_F\) requires access to all \(n^2\) entries.  
Similarly, checking whether the approximation satisfies
\[
H^{-1} A \approx I_n
\]
is not immediately efficient: while \(H^{-1}\) can be applied to a vector in \(O(n\log n)\) time using the hierarchical inverse, verifying the identity itself would still require \(O(n^2\log n)\) mat–vec applications.  
This remains, however, a significant improvement over the \(O(n^3)\) cost of dense inversion.

\begin{figure}[h]
    \centering 
    \begin{minipage}[b]{0.48\textwidth}
        \centering
        \includegraphics[width=\textwidth]{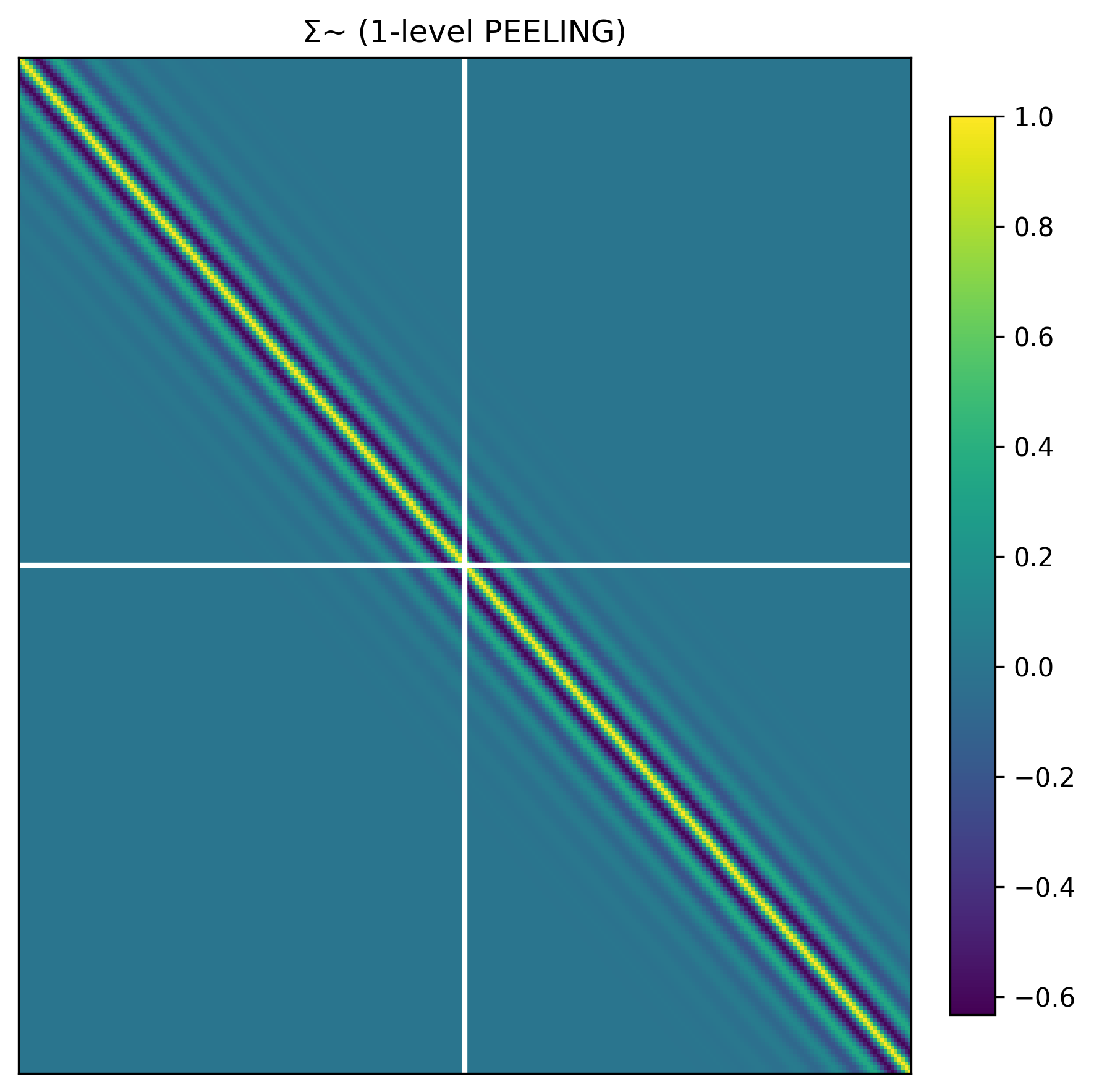}
        \caption{One-level HODLR approximation via the peeling algorithm. 
        The matrix is permuted so that the two leaf blocks are contiguous, exposing the coarse low-rank 
        left–right interaction captured at the first hierarchical level. 
        The one-level factorization attains a relative Frobenius error of \(6.95\times10^{-3}\),
        with eigenvalues of \(B=\tilde{\Sigma}^{-1}\Sigma\) spanning 
        \([0.6087,\,1.1520]\) and proxy KL divergence 
        \(D_{\mathrm{proxy}}=1.82\times10^{-2}\).}
        \label{fig:hodlr-peel-lvl1}
    \end{minipage}
    \hfill
    \begin{minipage}[b]{0.48\textwidth}
        \centering
        \includegraphics[width=\textwidth]{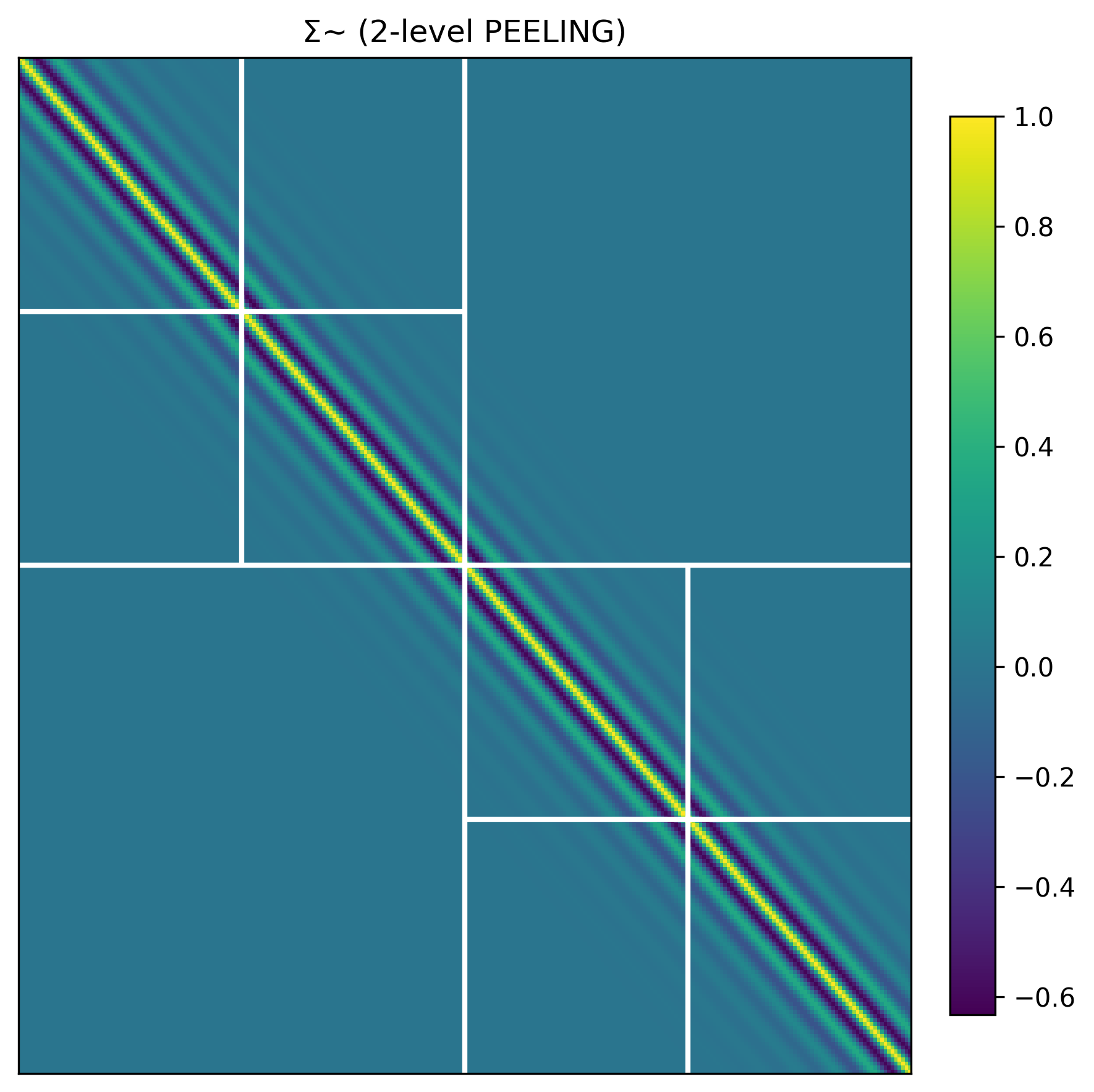}
        \caption{Two-level HODLR approximation, obtained by recursively 
        applying the peeling step to each child block. The four-leaf structure captures finer 
        off-diagonal dependencies while preserving hierarchical compression. 
        The two-level factorization yields a relative Frobenius error of \(9.26\times10^{-3}\),
        with eigenvalues of \(B=\tilde{\Sigma}^{-1}\Sigma\) spanning 
        \([0.6087,\,1.1773]\) and proxy KL divergence 
        \(D_{\mathrm{proxy}}=1.85\times10^{-2}\).}
        \label{fig:hodlr-peel-lvl2}
    \end{minipage}
\end{figure}

We combine the above HODLR construction with the subblock stochastic Lanczos quadrature (Subblock SLQ) framework to estimate divergence measures such as the Kullback–Leibler and Wasserstein–2 distances efficiently.  
Given a matrix \(A\) and its HODLR approximation \(\tilde A\) obtained via peeling, we form
\[
B = \tilde A^{-1} A,
\qquad
f(\lambda) = \lambda - \log \lambda - 1,
\]
and evaluate
\[
D_{\mathrm{proxy}} = \frac{1}{2}\operatorname{tr}\!\big(f(B)\big)
\]
using random principal subblocks of \(B\) as described in Algorithm~3.1.  
Since both \(A\) and \(\tilde A^{-1}\) admit fast hierarchical mat–vecs, this computation scales as \(O(n\log n)\) per iteration.  
In practice, we observe that the resulting proxy KL divergence closely matches the full-matrix divergence even for oscillatory kernels such as
\[
K(x_i, x_j) = e^{-\frac{|x_i-x_j|}{\ell}} \cos(2\pi \nu |x_i-x_j|),
\]
as illustrated in Figure~\ref{fig:hodlr-peel-lvl1} and~\ref{fig:hodlr-peel-lvl2}, demonstrating that proxy KL provides an effective low-memory surrogate for large-scale divergence estimation.

The accompanying implementation constructs one- and two-level HODLR factorizations of oscillatory exponential kernels using randomized SVD and verifies the approximation through proxy KL divergence. 

As shown in Figures~\ref{fig:hodlr-peel-lvl1}--\ref{fig:hodlr-peel-lvl2}, 
the peeling-based HODLR construction efficiently recovers both coarse and fine
low-rank off-diagonal structures of the oscillatory kernel using only 
\(O(n\log n)\) matrix–vector products.  
The proxy Kullback–Leibler divergence is evaluated using principal subblocks of 
\(B = \tilde{\Sigma}^{-1}\Sigma\), where each subblock \(B_s\) captures the local
spectral behavior of the hierarchical approximation.
In the experiments with \(n=256\), the one- and two-level HODLR factorizations
achieve relative Frobenius errors of order \(10^{-2}\) and proxy KL values near
\(1.8\times10^{-2}\), demonstrating strong agreement between hierarchical
and dense representations.

\begin{figure}[h]
    \centering
    \includegraphics[width=0.6\textwidth]{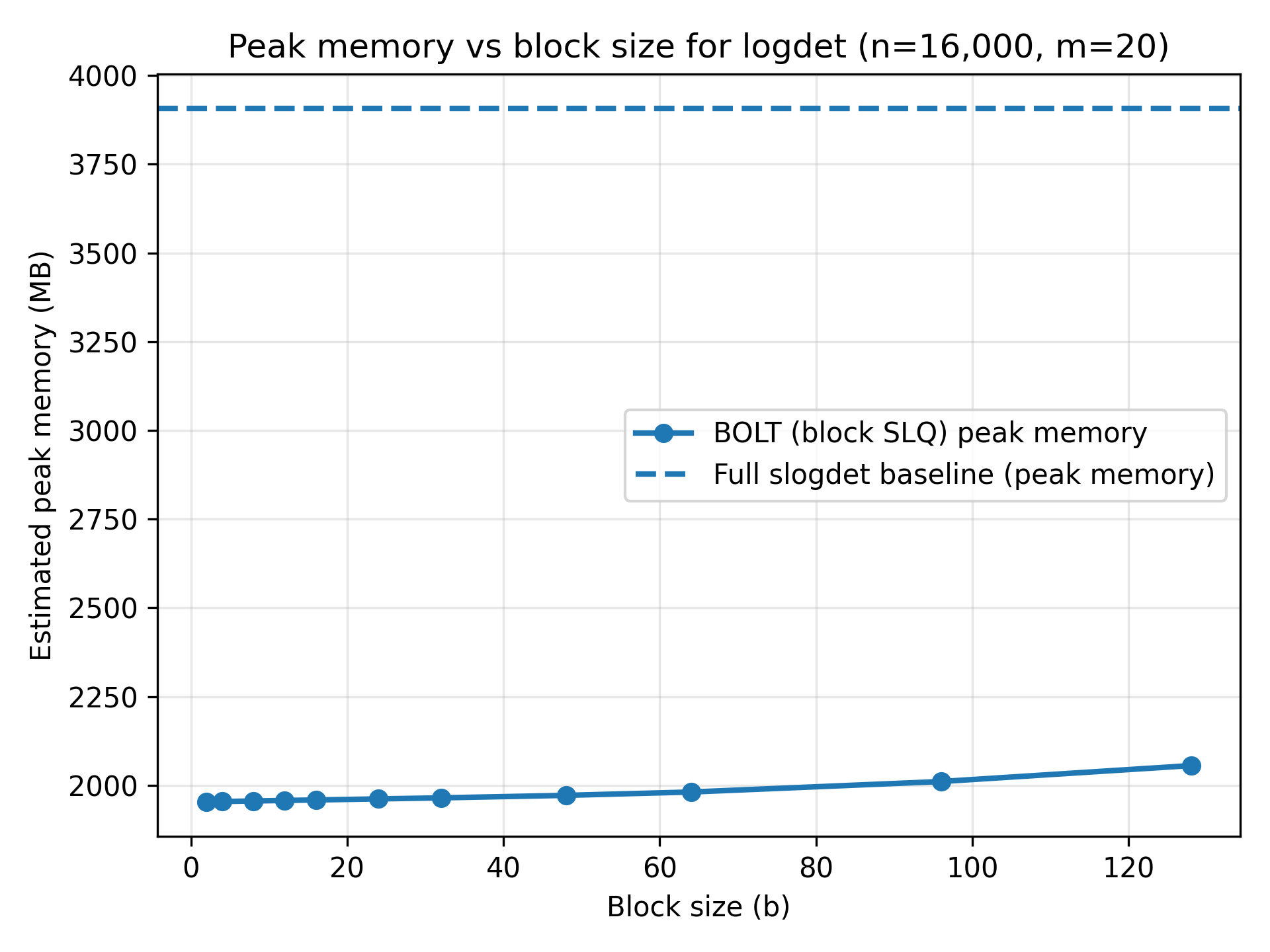}
    \caption{
    Peak memory usage of the BOLT estimator and the dense \texttt{slogdet} baseline for a kernel matrix of size \(n = 16{,}000\).
    BOLT achieves a nearly constant 2~GB footprint across block sizes, while the dense baseline requires about 4~GB.
    }
    \label{fig:memory_plot}
\end{figure}

To assess memory scalability, we applied the same experimental setup to a kernel matrix of size \(n = 16{,}000\).
Each BOLT experiment used twenty Lanczos iterations with block sizes \(b \in \{2, 4, \ldots, 128\}\), and the resulting peak memory usage is shown in Figure~\ref{fig:memory_plot}.
In the dense baseline, the Kullback--Leibler divergence between Gaussian models
\(\mathcal{N}(0, \Sigma)\) and \(\mathcal{N}(0, \widetilde{\Sigma})\) is computed via
\[
    D_{\mathrm{KL}}(\Sigma \,\|\, \widetilde{\Sigma})
    = \tfrac{1}{2}\bigl[
    \operatorname{tr}(\widetilde{\Sigma}^{-1}\Sigma)
    - n
    + \log\!\det(\widetilde{\Sigma})
    - \log\!\det(\Sigma)
    \bigr],
\]
where both log-determinant terms are evaluated using the dense \texttt{slogdet} factorization.
This procedure requires storing multiple \(n\times n\) arrays (\(\Sigma\), \(\widetilde{\Sigma}\), and intermediate factors) simultaneously, resulting in a peak memory footprint of roughly 3.9--4.0~GB for \(n=16{,}000\).
In contrast, the BOLT method maintains a nearly constant 2~GB requirement across block sizes, scaling linearly with the number of matrix--vector products while avoiding the dense factorization cost associated with $\log\!\det$.


\section{Conclusion}

We introduced BOLT, a block-orthonormal SLQ method, and Subblock SLQ, a trace estimator tailored for memory-limited settings with partial matrix access. Both reformulate KL divergence as a trace computation and offer scalable alternatives to classical estimators. BOLT matches Hutch++ in complexity but has a simpler implementation and avoids spectral decay assumptions, improving performance in flat-spectrum regimes. Subblock SLQ operates on small submatrices, making it ideal for large-scale problems with strict memory constraints. The same framework also supports Wasserstein-2 ($W_2$) distance estimation, which avoids matrix inversion and remains valid for positive semidefinite covariances.

This family of estimators is broadly applicable across covariance estimation~\cite{schafer2021sparse}, stochastic trace estimation~\cite{fuentes2023estimating, ubaru2017fast, kuczynski1992estimating, epperly2024xtrace}, and large-scale scientific computing~\cite{bai1996some}, where efficient KL evaluation is critical and full matrix access is often impractical. To our knowledge, this is the first unified KL framework based on trace functionals that handles full-rank, low-rank, and sparse regimes.

\section*{Acknowledgment}
 This work was supported by the U.S. Department of Energy, Office of Science (SC), Advanced Scientific Computing Research (ASCR), Competitive Portfolios Project on Energy Efficient Computing: A Holistic Methodology, under Contract DE-AC02-06CH11357.



\bibliographystyle{siamplain}
\bibliography{references}

\newpage
\appendix


\section{Proofs of Theoretical Results} \label{appendix:proofs}

\subsection{Proof of Theorem~\ref{thm:blockSLQ}} \label{appendix:blockSLQ}
\begin{proof}
We first prove unbiasedness. Since the columns of \(V\) are orthonormal and drawn from the Haar measure on the Stiefel manifold, we have
\[
\mathbb{E}[VV^T] = \frac{b}{n}I_n.
\]
Therefore,
\[
\mathbb{E}\Bigl[\operatorname{tr}\bigl(V^T f(A)V\bigr)\Bigr]
=\operatorname{tr}\Bigl(f(A)\,\mathbb{E}[VV^T]\Bigr)
=\frac{b}{n}\operatorname{tr}\bigl(f(A)\bigr).
\]
Defining
\[
X(V)=\frac{n}{b}\operatorname{tr}\bigl(V^T f(A)V\bigr),
\]
it follows immediately that
\[
\mathbb{E}[X(V)] = \operatorname{tr}\bigl(f(A)\bigr).
\]
In the extreme case \(b=n\), any orthonormal \(V\) is orthogonal (i.e. \(VV^T=I_n\)) so that
\[
\operatorname{tr}\bigl(V^T f(A)V\bigr)=\operatorname{tr}\bigl(f(A)\bigr),
\]
and hence \(X(V)\) is deterministic with zero variance.

Next, we derive the variance. Write the columns of \(V\) as \(v_1, v_2, \dots, v_b\), so that
\[
\operatorname{tr}\Bigl(V^T f(A)V\Bigr)=\sum_{r=1}^{b} v_r^T f(A)v_r.
\]
For a single scalar unit vector \(v\) uniformly distributed on the sphere in \(\mathbb{R}^n\), we express \(v\) in the eigenbasis of \(A\) by writing \(A=Q\Lambda Q^T\) with \(\Lambda = \mathrm{diag}(\lambda_1,\dots,\lambda_n)\) and \(v = Qz\). Then,
\[
v^T f(A)v = z^T f(\Lambda) z = \sum_{i=1}^n f(\lambda_i) z_i^2.
\]
Because \(z\) is uniformly distributed on the unit sphere, the moments of its components are \cite{mardia2009directional}:
\[
\mathbb{E}[z_i^2] = \frac{1}{n},\quad \mathbb{E}[z_i^4] = \frac{3}{n(n+2)},\quad \mathbb{E}[z_i^2 z_j^2] = \frac{1}{n(n+2)} \quad (i\neq j).
\]
Thus,
\[
\begin{aligned}
\mathbb{E}\Bigl[\Bigl(v^T f(A)v\Bigr)^2\Bigr]
&=\sum_{i=1}^n f(\lambda_i)^2\,\mathbb{E}[z_i^4]
+\sum_{i\neq j} f(\lambda_i)f(\lambda_j)\,\mathbb{E}[z_i^2z_j^2]\\[1mm]
&=\frac{3}{n(n+2)}\sum_{i=1}^n f(\lambda_i)^2
+\frac{1}{n(n+2)}\left[\left(\sum_{i=1}^n f(\lambda_i)\right)^2-\sum_{i=1}^n f(\lambda_i)^2\right]\\[1mm]
&=\frac{2}{n(n+2)}\sum_{i=1}^n f(\lambda_i)^2
+\frac{1}{n(n+2)}\left(\sum_{i=1}^n f(\lambda_i)\right)^2.
\end{aligned}
\]

Therefore, the variance of a scalar quadratic form is

\[
\begin{aligned}
\operatorname{Var}\bigl[v^T f(A)v\bigr]
&=\mathbb{E}\Bigl[\bigl(v^T f(A)v\bigr)^2\Bigr]
-\Bigl(\frac{1}{n}\sum_{i=1}^n f(\lambda_i)\Bigr)^2\\[3pt]
&=\frac{2}{n(n+2)}\sum_{i=1}^n f(\lambda_i)^2
-\frac{2}{n^2(n+2)}\Bigl(\sum_{i=1}^n f(\lambda_i)\Bigr)^2.
\end{aligned}
\]

Now, consider the block estimator. Define
\[
S = \sum_{r=1}^b v_r^T f(A)v_r.
\]
If the \(b\) vectors were independent, we would have 
\[
\operatorname{Var}(S)= b\,\operatorname{Var}\bigl[v^T f(A)v\bigr].
\]
However, because the columns of \(V\) are orthonormal, for any \(i \neq j\) we have
\[
\operatorname{Cov}\Bigl(v_i^T f(A)v_i,\, v_j^T f(A)v_j\Bigr)
=-\frac{1}{n-1}\,\operatorname{Var}\Bigl(v^T f(A)v\Bigr).
\]
Thus,
\[
\begin{aligned}
\operatorname{Var}(S)
&=\sum_{r=1}^b \operatorname{Var}\bigl[v^T f(A)v\bigr] + \sum_{r\neq s}\operatorname{Cov}\Bigl(v_r^T f(A)v_r,\, v_s^T f(A)v_s\Bigr)\\[1mm]
&= b\,\operatorname{Var}\bigl[v^T f(A)v\bigr] - \frac{b(b-1)}{n-1}\,\operatorname{Var}\bigl[v^T f(A)v\bigr]\\[1mm]
&=\operatorname{Var}\bigl[v^T f(A)v\bigr]\left(b-\frac{b(b-1)}{n-1}\right).
\end{aligned}
\]
We call this the {\em self-averaging effect}. Since 
\[
X(V)=\frac{n}{b}S,
\]
we obtain
\[
\operatorname{Var}\bigl[X(V)\bigr] = \left(\frac{n}{b}\right)^2 \operatorname{Var}(S)
=\left(\frac{n}{b}\right)^2 \operatorname{Var}\bigl[v^T f(A)v\bigr]\left(b-\frac{b(b-1)}{n-1}\right).
\]
Substituting the expression for \(\operatorname{Var}\bigl[v^T f(A)v\bigr]\) yields
\[
\operatorname{Var}\bigl[X(V)\bigr] = \left(\frac{n}{b}\right)^2 \frac{2}{n(n+2)}\left(\sum_{i=1}^n f(\lambda_i)^2-\frac{1}{n}\Bigl(\operatorname{tr}(f(A))\Bigr)^2\right)\left(b-\frac{b(b-1)}{n-1}\right).
\]
Simplifying, we note that
\[
\left(\frac{n}{b}\right)^2 \frac{2}{n(n+2)} = \frac{2n}{b^2(n+2)}.
\]
Thus, 
\[
\operatorname{Var}\bigl[X(V)\bigr] = \frac{2n}{b^2(n+2)}\left(b-\frac{b(b-1)}{n-1}\right)
\left(\sum_{i=1}^n f(\lambda_i)^2-\frac{1}{n}\Bigl(\operatorname{tr}(f(A))\Bigr)^2\right).
\]
We can factor a \(b\) from the term in parentheses:
\[
b-\frac{b(b-1)}{n-1} = b\left(1-\frac{b-1}{n-1}\right),
\]
so that
\[
\operatorname{Var}\bigl[X(V)\bigr] = \frac{2n}{b^2(n+2)}\,b\left(1-\frac{b-1}{n-1}\right)
\left(\sum_{i=1}^n f(\lambda_i)^2-\frac{1}{n}\Bigl(\operatorname{tr}(f(A))\Bigr)^2\right).
\]
Cancel one \(b\) to obtain
\[
\operatorname{Var}\bigl[X(V)\bigr] = \frac{2n}{b(n+2)}\left(1-\frac{b-1}{n-1}\right)
\left(\sum_{i=1}^n f(\lambda_i)^2-\frac{1}{n}\Bigl(\operatorname{tr}(f(A))\Bigr)^2\right). 
\]

Where, 
\[
\frac{2n}{b(n+2)}\Bigl(1-\frac{b-1}{n-1}\Bigr)
\;=\;
-\frac{2n}{(n+2)(n-1)}
\;+\;\frac{2n^2}{(n+2)(n-1)}\,\frac{1}{b}
\;+\;0\cdot\frac{1}{b^2}
\;+\;\mathcal{O}\!\bigl(b^{-3}\bigr)
\;=\;
\mathcal{O}\bigl(b^{-1}\bigr)
\]

Averaging over \(q\) independent probes gives
\[
\operatorname{Var}[\widehat{T}] = \frac{1}{q}\operatorname{Var}\bigl[X(V)\bigr].
\]

Therefore the variance decays like $\mathcal{O}\bigl(\frac{1}{qb}\bigr)$.


Now we derive a lower bound on the variance of Hutch++ with Gaussian probes for comparison. Set \(k=\lfloor m/3\rfloor\), and assume for the purpose of a worst-case lower bound that the sketch matrix \(Q \in \mathbb{R}^{n \times k}\) captures exactly the top-\(k\) eigenspace of \(A\). Then the residual is
\[
R = (I - QQ^T) A (I - QQ^T) = A - A_k,
\]
where \(A_k\) is the best rank-\(k\) approximation to \(A\) in the Frobenius norm. Under this assumption, the Hutch++ estimator decomposes as
\[
\widehat T_{\mathrm{H++}}(A) = \underbrace{\operatorname{tr}(Q^T A Q)}_{\text{deterministic}} + \frac{3}{m} \sum_{i=1}^{m/3} g_i^T R g_i,
\]
where the vectors \(g_i \sim \mathcal{N}(0, I_n)\) are independent. Since the first term is deterministic, the variance comes entirely from the residual component:
\[
\operatorname{Var}[\widehat T_{\mathrm{H++}}(A)] = \left( \frac{3}{m} \right)^2 \cdot \frac{m}{3} \cdot \operatorname{Var}[g^T R g] = \frac{3}{m} \operatorname{Var}[g^T R g].
\]
It is well known that for a Gaussian quadratic form, \(\operatorname{Var}[g^T R g] = 2 \|R\|_F^2\) (see \cite{meyer2021hutch++}), and thus
\[
\operatorname{Var}[\widehat T_{\mathrm{H++}}(A)] \ge \frac{6}{m} \|A - A_k\|_F^2 = \frac{6}{m}\sum_{i=k+1}^n\lambda_i^2.
\]

By contrast, the variance of our BOLT estimator satisfies
\[
\operatorname{Var}[X(V)] = \frac{2n}{b(n+2)}\left(1 - \frac{b-1}{n-1} \right) \left( \sum_{i=1}^n \lambda_i^2 - \frac{1}{n} \left( \sum_{i=1}^n \lambda_i \right)^2 \right).
\]
Taking the ratio of BOLT's variance to Hutch++’s worst-case lower bound gives
\begin{equation}
\label{eq:bolt-vs-hpp-var}
\frac{\operatorname{Var}[X(V)]}{\operatorname{Var}[\widehat T_{\mathrm{H++}}(A)]}
\le
\frac{mn}{3b(n+2)}
\left(1 - \frac{b-1}{n-1} \right)
\cdot
\frac{\sum_{i=1}^n \lambda_i^2 - \frac{1}{n} \left( \sum_{i=1}^n \lambda_i \right)^2}
     {\sum_{i=k+1}^n \lambda_i^2}.
\end{equation}
When the spectrum is flat, nearly flat, or slowly decaying, so that the residual \(\|A - A_k\|_F^2\) remains large, BOLT can significantly outperform Hutch++ in variance. Conversely, when the spectrum exhibits rapid decay and the dominant spectral mass is concentrated in the top \(k\) eigenvalues, Hutch++ achieves lower variance by effectively projecting out the high-energy components.


Finally, since each \(X(V)\) is almost surely contained in \([\ell_{\min},\ell_{\max}]\), Hoeffding's inequality implies that for any \(\varepsilon>0\)
\[
\mathbb{P}\!\Biggl(\Bigl|\widehat{T}-\operatorname{tr}(f(A))\Bigr|\ge\varepsilon\Biggr)
\le 2\exp\!\Biggl(-\frac{2q\,\varepsilon^2}{(\ell_{\max}-\ell_{\min})^2}\Biggr).
\]
Equivalently, with probability at least \(1-\delta\),
\[
\Bigl|\widehat{T}-\operatorname{tr}(f(A))\Bigr|\le (\ell_{\max}-\ell_{\min})\sqrt{\frac{\ln(2/\delta)}{2q}}.
\]
\end{proof}

\subsection{Proof of Lemma~\ref{thm:fullrank}} \label{appendix:fullrank}
\begin{proof}
Let \(P_S\in\{0,1\}^{n\times s}\) be the coordinate‐selection matrix that extracts the \(s\) indices in \(S\).  For each Gaussian sample \(u_i\sim\mathcal N(0,\Sigma)\), the projected vector
\[
v_i \;=\; P_S^T u_i
\]
is an \(s\)-dimensional Gaussian.  Its mean is 
\(\mathbb{E}[v_i] = P_S^T\,\mathbb{E}[u_i] = 0\), 
and its covariance is
\[
\mathbb{E}[v_i v_i^T]
= P_S^T\,\mathbb{E}[u_i u_i^T]\,P_S
= P_S^T\,\Sigma\,P_S.
\]
Hence 
\[
v_1,\dots,v_m \sim \text{iid}\;\mathcal N\bigl(0,\;P_S^T\Sigma P_S\bigr).
\]
By definition, the unnormalized sample covariance of these projections is
\[
\widetilde\Sigma_S
= \sum_{i=1}^m v_i v_i^T
= \sum_{i=1}^m (P_S^T u_i)(P_S^T u_i)^T
= P_S^T\!\Bigl(\sum_{i=1}^m u_i u_i^T\Bigr)P_S
= P_S^T\,\widetilde\Sigma \,P_S,
\]
which is a central Wishart matrix
\(\;W_s\bigl(P_S^T\Sigma P_S,\;m\bigr)\).

A key property of the Wishart distribution is that if the degrees of freedom \(m\) satisfy \(m \ge s\) and the scale matrix \(P_S^T\Sigma P_S\) is full rank (which it is, since \(\Sigma\succ0\)), then \(\widetilde\Sigma_S\) is invertible with probability one.  Equivalently,
\[
\Pr\bigl(\det\widetilde\Sigma_S=0\bigr)=0,
\]
so \(\text{rank}(\widetilde\Sigma_S)=s\) almost surely.  Conversely, if \(m<s\), then \(\widetilde\Sigma_S\) has rank at most \(m\), and so the transition from full‐rank to singular behavior occurs precisely at \(s=m\).
\end{proof}

\subsection{Proof of Lemma~\ref{lem:unbiasedsubblockbound}} \label{appendix:unbiasedsubblockbound}
\begin{proof}
First, for any \(i\in\{1,\dots,n\}\), \(\Pr(i\in S)=s/n\), so
\[
\mathbb{E}[\tr(f(A)_S)]
=\sum_{i=1}^n \Pr(i\in S)\,f(A)_{ii}
=\tfrac{s}{n}\,\tr(f(A)),
\]
and hence \(\mathbb{E}[X(S)]=\tr(f(A))\).


\end{proof}

\subsection{Proof of Theorem~\ref{thm:localization}}\label{appendix:localization} 
\begin{proof}
Write the Chebyshev expansion
\[
p_m(A)
=\sum_{k=0}^m c_k\,A^k.
\]
We will show that for each \(0\le k\le m\) and every \(i,j\in S\),
\[
[A^k]_{ij}
\;=\;
\bigl[(A_{S_r,S_r})^k\bigr]_{ij}.
\]
Indeed, by the definition of matrix powers,
\[
[A^k]_{ij}
=
\sum_{t_1,\dots,t_{k-1}=1}^n
A_{i,t_1}\,A_{t_1,t_2}\,\cdots\,A_{t_{k-1},j}.
\]
Each nonzero term in this sum corresponds to a walk of length \(k\) in \(G\) from \(i\) to \(j\) through vertices \(t_1,\dots,t_{k-1}\).  Since \(i\in S\) and \(k\le r\), every intermediate vertex and the endpoint \(j\) lie within distance \(\le k\le r\) of \(S\), hence belong to \(S_r\).  Thus all indices in any nonzero product lie in \(S_r\), and we may restrict each summation to \(t_\ell\in S_r\).  But that restricted sum is exactly the \((i,j)\)-entry of \((A_{S_r,S_r})^k\).  Therefore
\[
[A^k]_{ij}
=\bigl[(A_{S_r,S_r})^k\bigr]_{ij},
\]
as claimed.

Summing these equalities over \(k=0,\dots,m\) with weights \(c_k\) gives
\[
\bigl[p_m(A)\bigr]_{S,S}
=\sum_{k=0}^m c_k\,[A^k]_{S,S}
=\sum_{k=0}^m c_k\bigl[(A_{S_r,S_r})^k\bigr]_{S,S}
=\bigl[p_m(A_{S_r,S_r})\bigr]_{S,S}.
\]
Finally, since \(\deg p_m\le r\), no walk of length \(\le m\le r\) can leave \(S\) once it has returned; restricting further from the \(S_r\)-block to the \(S\)-block recovers
\(\bigl[p_m(A_{S_r,S_r})\bigr]_{S,S}=p_m(A_S)\).

To bound \(\Delta_S\), insert and subtract \(p_m\) twice:
\[
\Delta_S
=\bigl[\tr(f(A)_S)-\tr(p_m(A)_S)\bigr]
+\bigl[\tr(p_m(A)_S)-\tr(p_m(A_S))\bigr]
+\bigl[\tr(p_m(A_S))-\tr(f(A_S))\bigr].
\]
The middle term vanishes by the block-wise equality just shown.  Meanwhile the uniform approximation \(\|f-p_m\|_\infty\le\varepsilon\) implies
\(\|f(A)-p_m(A)\|_2\le\varepsilon\) and similarly on the submatrix.  Hence each of the two endpoint differences in trace is bounded by \(s\,\varepsilon\), giving
\(\lvert\Delta_S\rvert\le s\,\varepsilon +0+s\,\varepsilon=2s\,\varepsilon\).
\end{proof}


\section{FLOPS count of various matrix trace estimator}
\label{appendix:trace_flops}

The first three estimators, XTrace, XNysTrace~\cite{epperly2024xtrace}, and Hutch++, all build on the idea of combining low-rank sketching with randomized probing, but differ in how they form the sketch and treat the residual. These structural differences directly influence the computational cost of each method, as summarized in Table~\ref{table:flop-old}.

\begin{table}[h]
\centering
\caption{Flop counts to achieve \(\varepsilon\)-accuracy (\(m, b \sim \Theta(\varepsilon^{-1})\)) (\textbf{Na\"ive}).}
\label{table:flop-old}
\begin{tabular}{l c c c}
\toprule
Method      & \# Mat–vecs & Overhead (post–mat–vec)           & Total flops                            \\
\midrule
XTrace (Na\"ive)     & \(m\)  
            & \(O(N\,m^3 + m^4 + N\,m^2)\)       & \(O(N^2 m + N\,m^3 + m^4 + N\,m^2)\)   \\

XNysTrace  (Na\"ive)  & \(m\)  
            & \(O(m^4 + N\,m)\)                 & \(O(N^2 m + m^4 + N\,m)\)             \\

Hutch++     & \(m\)  
            & \(O(N\,m^2)\)                     & \(O(N^2 m + N\,m^2)\)                 \\

BOLT        & \(m\)  
            & \(O(N\,m^2)\)                     & \(O(N^2 m + N\,m^2)\)                 \\
\bottomrule
\end{tabular}
\end{table}

XTrace enforces the exchangeability principle by constructing a leave-one-out low-rank approximation for each of the \(m/2\) probing directions. It draws \(\omega_1,\dots,\omega_{m/2}\) and forms the sketch matrix \(Y = A[\omega_1\,\cdots\,\omega_{m/2}]\in\mathbb{R}^{N\times(m/2)}\). For each \(i\), the \(i\)th column is removed to obtain \(Y_{-i}\), and an orthonormal basis \(Q^{(i)} = \mathrm{orth}(Y_{-i})\) is computed. The trace estimator is
\[
\hat t_i = \tr\bigl({Q^{(i)}}^T A Q^{(i)}\bigr) + \omega_i^T (I - Q^{(i)}{Q^{(i)}}^T) A (I - Q^{(i)}{Q^{(i)}}^T) \omega_i,
\]
and the final output is \(\widehat{\tr}(A) = \frac{2}{m} \sum_{i=1}^{m/2} \hat t_i\). Computing all QR factorizations costs \(O(N m^3)\), pseudoinverses contribute \(O(m^4)\), and the trace and residual evaluations require \(O(N m^2)\).

XNysTrace modifies this by using Nyström approximations. It constructs
\[
A\langle\Omega_{-i}\rangle = A \Omega_{-i} \left( \Omega_{-i}^T A \Omega_{-i} \right)^{\dagger} (\Omega_{-i}^T A)^T,
\quad
\hat t_i = \tr(A\langle\Omega_{-i}\rangle) + \omega_i^T (A - A\langle\Omega_{-i}\rangle) \omega_i,
\]
requiring an SVD or Cholesky on each \((m{-}1)\times(m{-}1)\) matrix, yielding \(O(m^4)\) total overhead. The rest of the cost mirrors that of XTrace.

\paragraph{Efficient Implementation via Rank-One Update.}
In Section~2.1 of~\cite{epperly2024xtrace}, a more efficient implementation of the XTrace estimator is developed. Let \(A \in \mathbb{R}^{N\times N}\) and define the test matrix \(\Omega = [\omega_1\,\cdots\,\omega_{m/2}] \in \mathbb{R}^{N \times (m/2)}\). One first forms the sketch \(Y = A\Omega\) and performs the QR decomposition \(Y = QR\). Then, following the analysis in Appendix A.2 of~\cite{epperly2024xtrace}, the leave-one-out orthogonal projectors are computed by
\[
Q^{(i)} Q^{(i)T} = Q(I - s_i s_i^T)Q^T,
\]
where each \(s_i\in \mathbb{R}^m\) is a unit vector in the nullspace of \(R^T_{-i}\). All vectors \(s_1,\dots,s_{m/2}\) can be constructed simultaneously via
\[
S = R^{-*} D,
\]
where \(D\) is a diagonal matrix used to normalize the columns of \(S\). This reduces the cost of all leave-one-out projectors to \(O(m^3)\) total.

\paragraph{Hutch++.}
Hutch++ achieves variance reduction without leave-one-out constructions by splitting the total of \(m\) mat–vecs into three groups of size \(s = m/3\). The first \(s\) probes are used to build a low-rank sketch \(Y = A[\omega_{s+1}, \dots, \omega_{2s}]\), and an orthonormal basis \(Q = \mathrm{orth}(Y)\) is computed. This basis is then used to exactly compute the projected trace \(\tr(Q^T A Q)\). The remaining \(s\) probes are used to estimate the residual trace \(\tr((I - QQ^T) A (I - QQ^T))\) via a Hutchinson-type average. The QR factorization of the \(N \times s\) matrix \(Y\) costs \(O(N s^2) = O(N m^2)\), and the subsequent evaluation of both the projected trace and residual adds another \(O(N m^2)\), since each involves applying \(A\) and performing inner products with matrices of size \(N \times s\). The total mat–vec cost is again \(O(N^2 m)\), as \(A\) is applied to all \(m\) probe vectors. Therefore, the total cost is \(O(N^2 m + N m^2)\), significantly lower than the leave-one-out variants due to the reuse of the same sketching basis across all evaluations.

\noindent
While all methods share a leading matrix–vector cost of \(O(N^2 m)\), XTrace and XNysTrace incur an additional \(O(m^3)\) overhead from the rank-one update implementation~\cite{epperly2024xtrace}. These extra costs make Hutch++ and BOLT more scalable, especially in settings where subblock extensions or memory constraints are relevant. A summary of the computational complexity under this optimized implementation is provided in Table~\ref{table:flop}.


\section{Eigenvalue Distributions for Wishart Matrices}
\label{appendix:eigenvalue_distributions}

We summarize key results about the eigenvalue distribution of Wishart matrices, which explain the rank transition behavior discussed in the main text.

Let \(\widetilde{\Sigma} = \sum_{i=1}^m u_i u_i^T\), where \(u_i \sim \mathcal{N}(0, \Sigma)\) and \(\Sigma \in \mathbb{R}^{n \times n}\) is full rank. Then \(\widetilde{\Sigma} \sim W_n(\Sigma, m)\) is a sample covariance matrix of rank at most \(m\). For any index set \(S\subset\{1,\dots,n\}\) of size \(s\), the subblock \(\widetilde{\Sigma}_S = P_S^T \widetilde{\Sigma} P_S \sim W_s(P_S^T \Sigma P_S, m)\) is itself Wishart-distributed. The smallest eigenvalue of \(\widetilde{\Sigma}_S\) undergoes a sharp phase transition at \(s = m\): subblocks are full rank almost surely for \(s \le m\), while for \(s > m\), they become singular with probability one. This transition governs when subblock-based KL estimation remains well-posed.

We recall the following classical result:

\begin{theorem}[Smallest Eigenvalue Distribution for Wishart Matrices with Gaussian Entries (Theorem~4.3, \cite{edelman1988eigenvalues})]
\label{thm:edelman}
Let \(W(m,n)\) be a Wishart matrix, formed from an \(m\times n\) matrix with independent, mean-zero Gaussian entries. Then the probability density function (pdf) of the smallest eigenvalue \(\lambda_{\min}\) is given by
\[
f(\lambda) = c_{m,n}\,\lambda^{\frac{n-m-1}{2}}\,e^{-\lambda/2}\,g(\lambda),
\]
where
\[
g(\lambda) =
\begin{cases}
P_{m,n}(\lambda), & \text{if } n-m \text{ is odd},\\
Q_{m,n}(\lambda)\, U\!\Bigl(\frac{m-1}{2},\frac{1}{2},\lambda\Bigr) + R_{m,n}(\lambda)\, U'\!\Bigl(\frac{m-1}{2},\frac{1}{2},\lambda\Bigr), & \text{if } n-m \text{ is even}.
\end{cases}
\]
Here, \(U(a,b,\lambda)\) denotes the Tricomi confluent hypergeometric function, and \(U'\) its derivative with respect to \(\lambda\).
The constant \(c_{m,n}\) depends explicitly on \(m\) and \(n\) and products of Gamma functions.
\end{theorem}

\begin{remark}
Although Theorem~\ref{thm:edelman} assumes Gaussian entries, similar eigenvalue statistics hold for a broad class of random matrices satisfying finite-moment assumptions, by universality results such as those of Tao and Vu~\cite{tao2010random}.
\end{remark}

In the special case where \(W(m,m)\) is a square Wishart matrix and \(m\to\infty\), the scaled smallest eigenvalue \(m\lambda_{\min}\) converges in distribution to a density
\[
f(x) = \frac{1+\sqrt{x}}{2\sqrt{x}}\,e^{-(x/2+\sqrt{x})}.
\]

\begin{figure}[h]
    \centering
    \begin{minipage}[b]{0.48\textwidth}
        \centering
        \includegraphics[width=\textwidth]{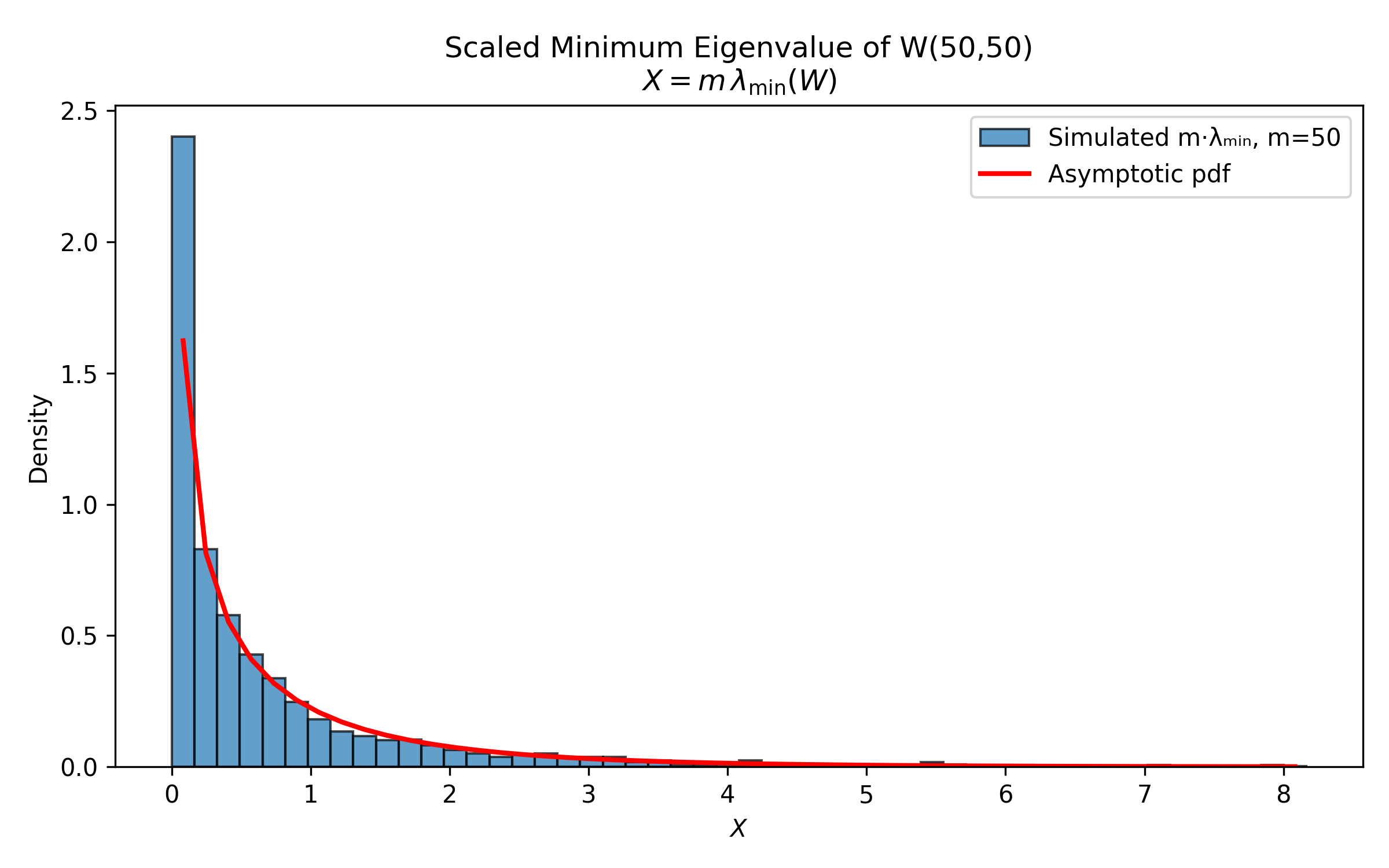}
        \caption{Histogram of scaled smallest eigenvalue of \(W(m,m)\) for \(m = 50\).}
        \label{fig:scaled_min_eig_hist}
    \end{minipage}
    \hfill
    \begin{minipage}[b]{0.48\textwidth}
        \centering
        \includegraphics[width=\textwidth]{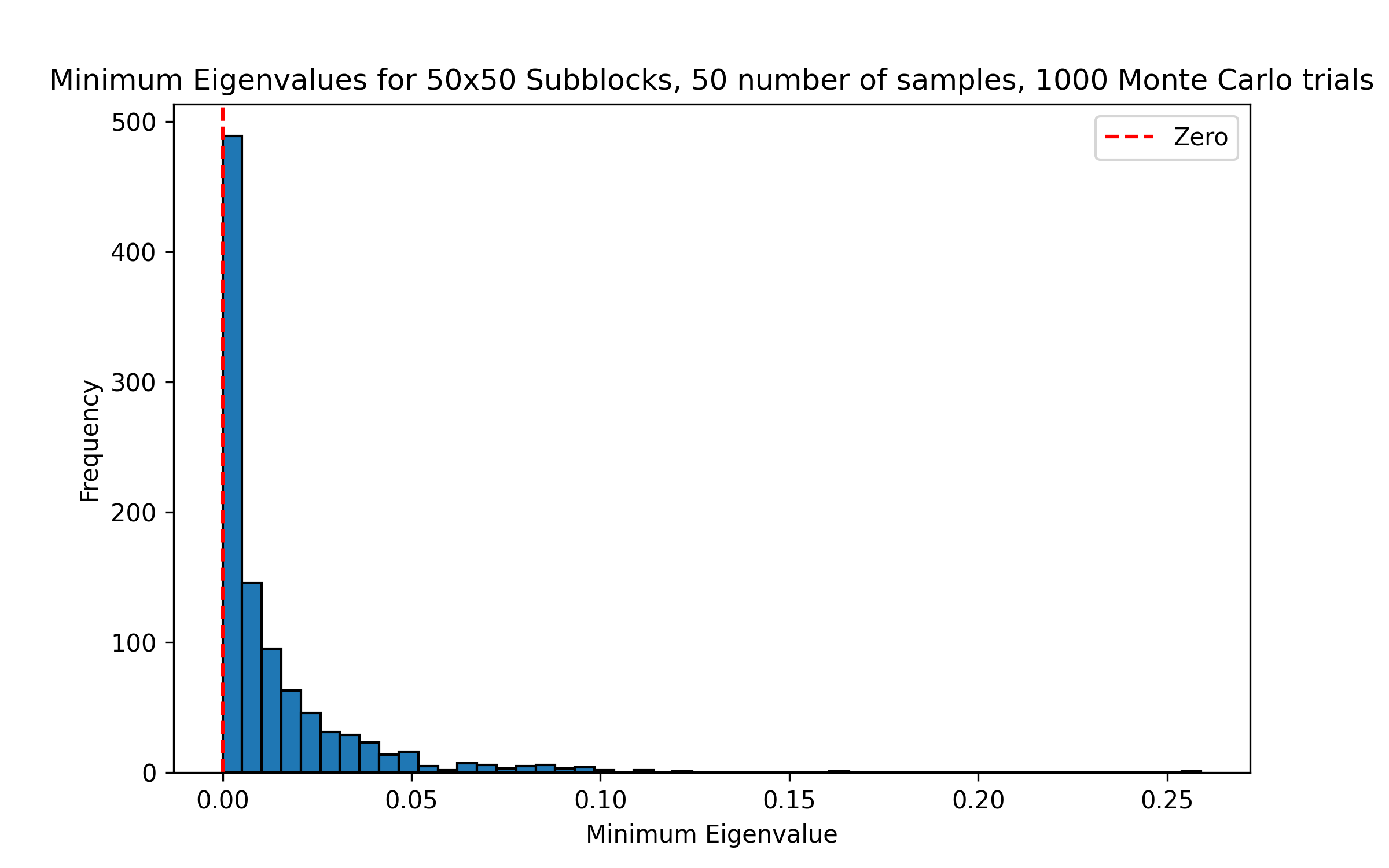}
        \caption{Histogram of smallest eigenvalue of \(A_s\) when \(s = m\).}
        \label{fig:hist1}
    \end{minipage}
\end{figure}

\begin{figure}[h]
    \centering
    \begin{minipage}{0.48\textwidth}
        \centering
        \includegraphics[width=\textwidth]{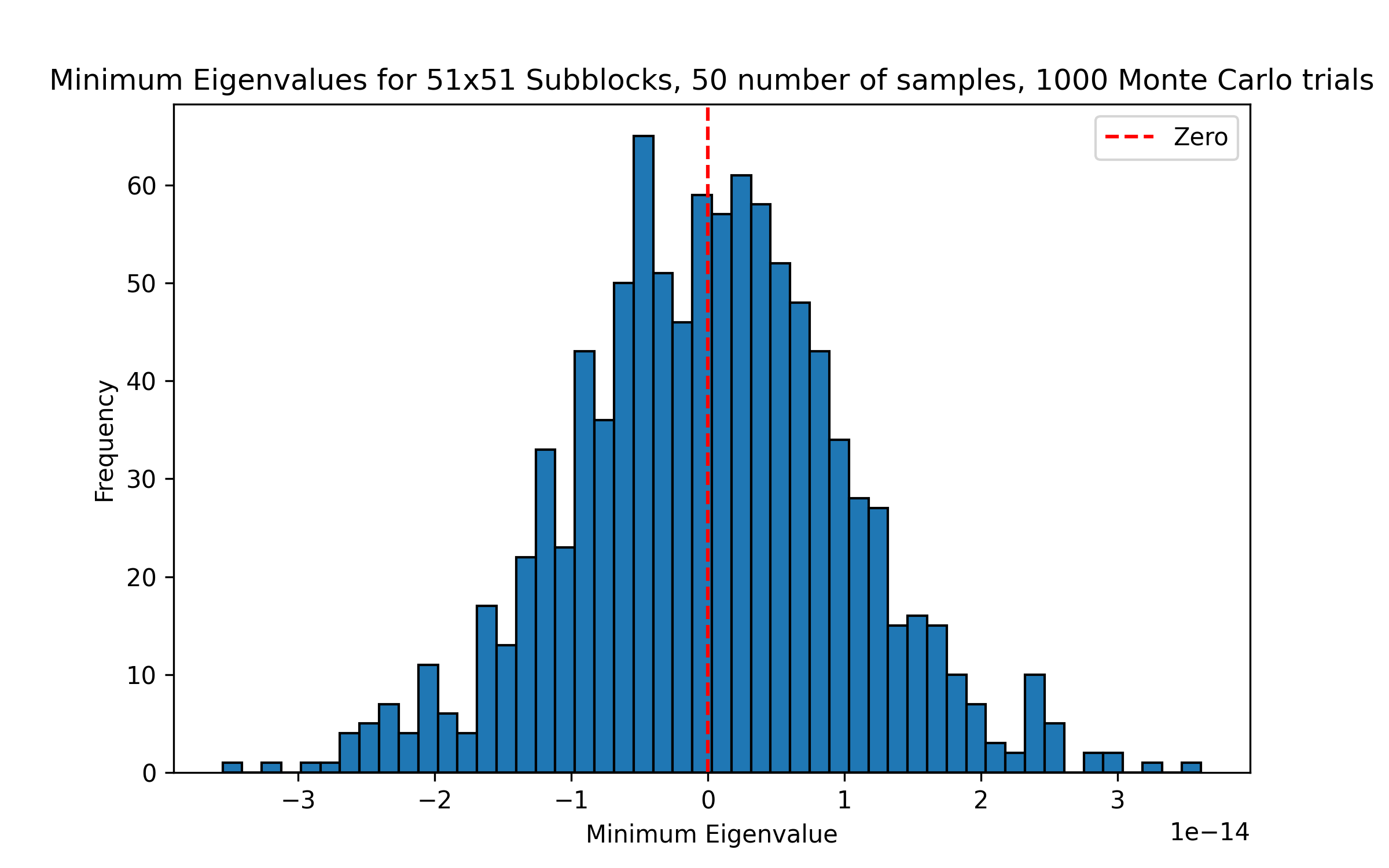}
        \caption{Histogram of the smallest eigenvalue of \(A_s\) when \(s>m\).}
        \label{fig:hist2}
    \end{minipage}
    \hfill
    \begin{minipage}{0.48\textwidth}
        \centering
        \includegraphics[width=\textwidth]{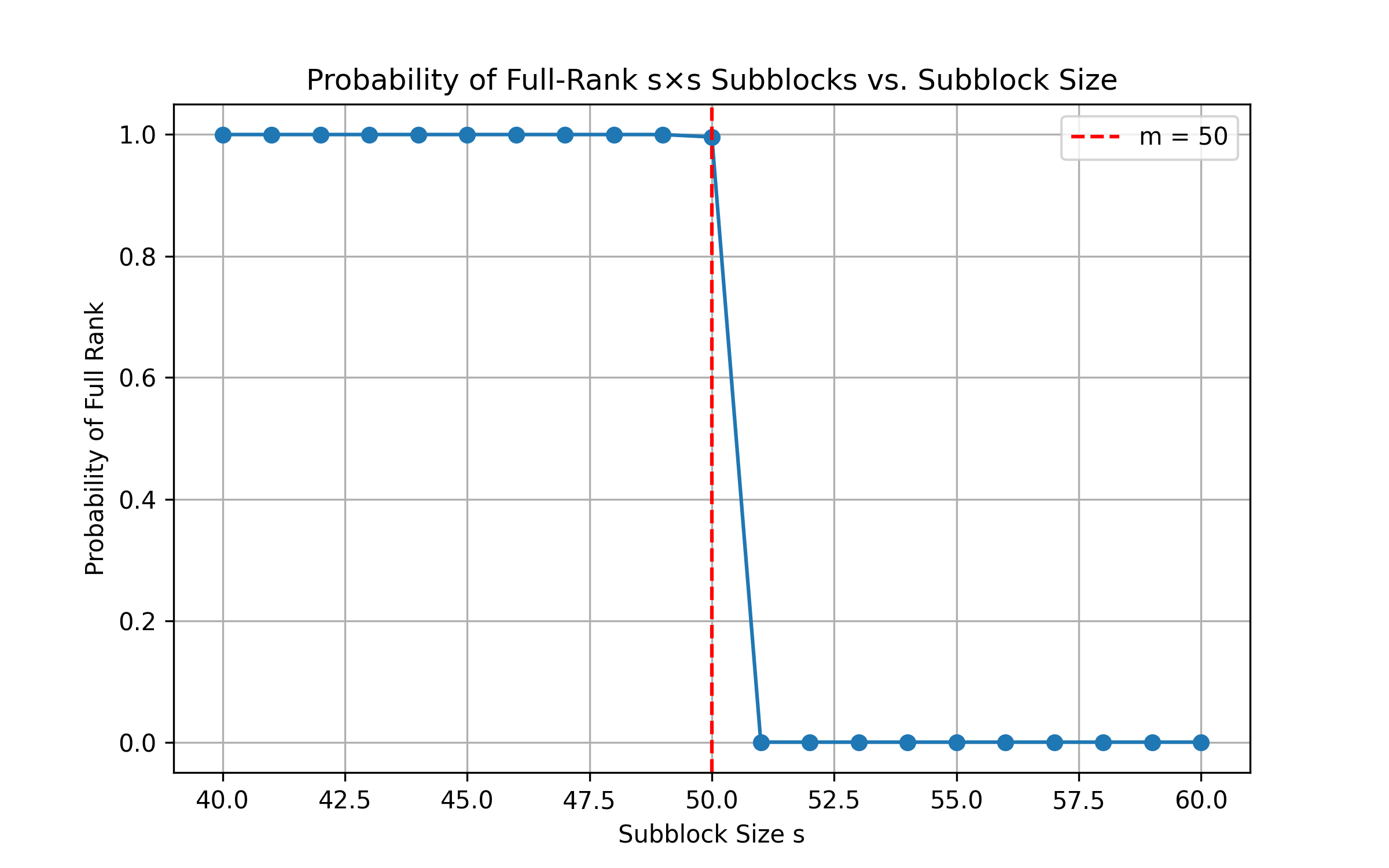}
        \caption{Probability of different subblock sizes being full rank.}
        \label{fig:prob_fullrank}
    \end{minipage}
\end{figure}

Figure~\ref{fig:scaled_min_eig_hist} shows the histogram of scaled smallest eigenvalues for \(m=50\), while Figures~\ref{fig:hist1}–\ref{fig:hist2} illustrate the transition from full-rank to singular subblocks as \(s\) exceeds \(m\).

\section{Numerical Experiment Details} \label{appendix:experiment-details}

\paragraph{Synthetic Cholesky Estimation.}
We generate a covariance matrix $\Sigma$ using a Gaussian kernel with bandwidth $\sigma = 2.0$ and simulate $m$ samples from $\mathcal{N}(0, \Sigma)$. For the full-rank case, we use $m = 10000$, and for the singular case, $m = 50$. Let $Y \in \mathbb{R}^{m \times n}$ denote the matrix of samples. We estimate the Cholesky factor $L$ by solving for $P = LL^T$ such that $P\Sigma \approx I$, using Newton's method on each column. When $m < n$, the estimate becomes low-rank, and standard KL divergence becomes undefined.

\begin{figure}[h]
    \centering
    \includegraphics[width=1\textwidth]{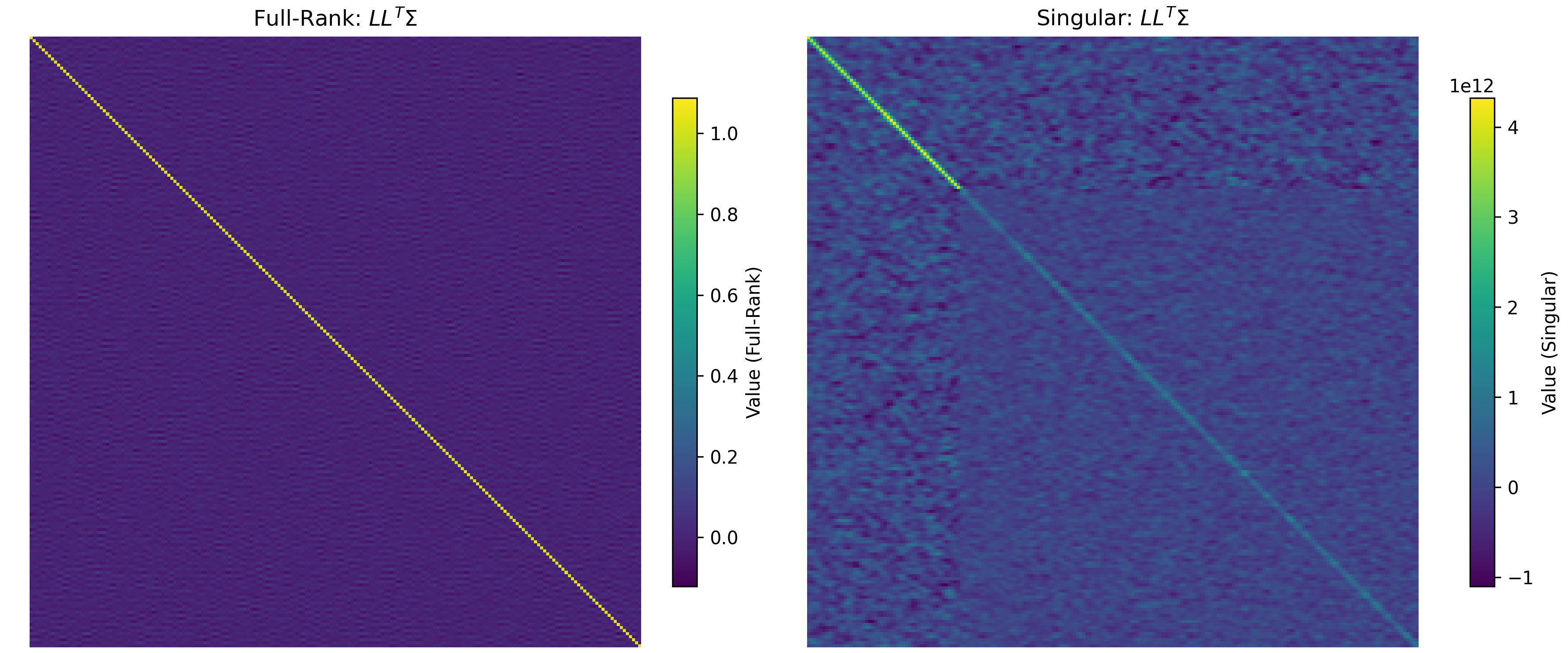}
    \caption{Full-Rank KL: 1.05 \\ Singular Case KL: 1.07e+14}
    \label{fig:err}
\end{figure}

\begin{figure}[h]
    \centering
    \includegraphics[width=1\textwidth]{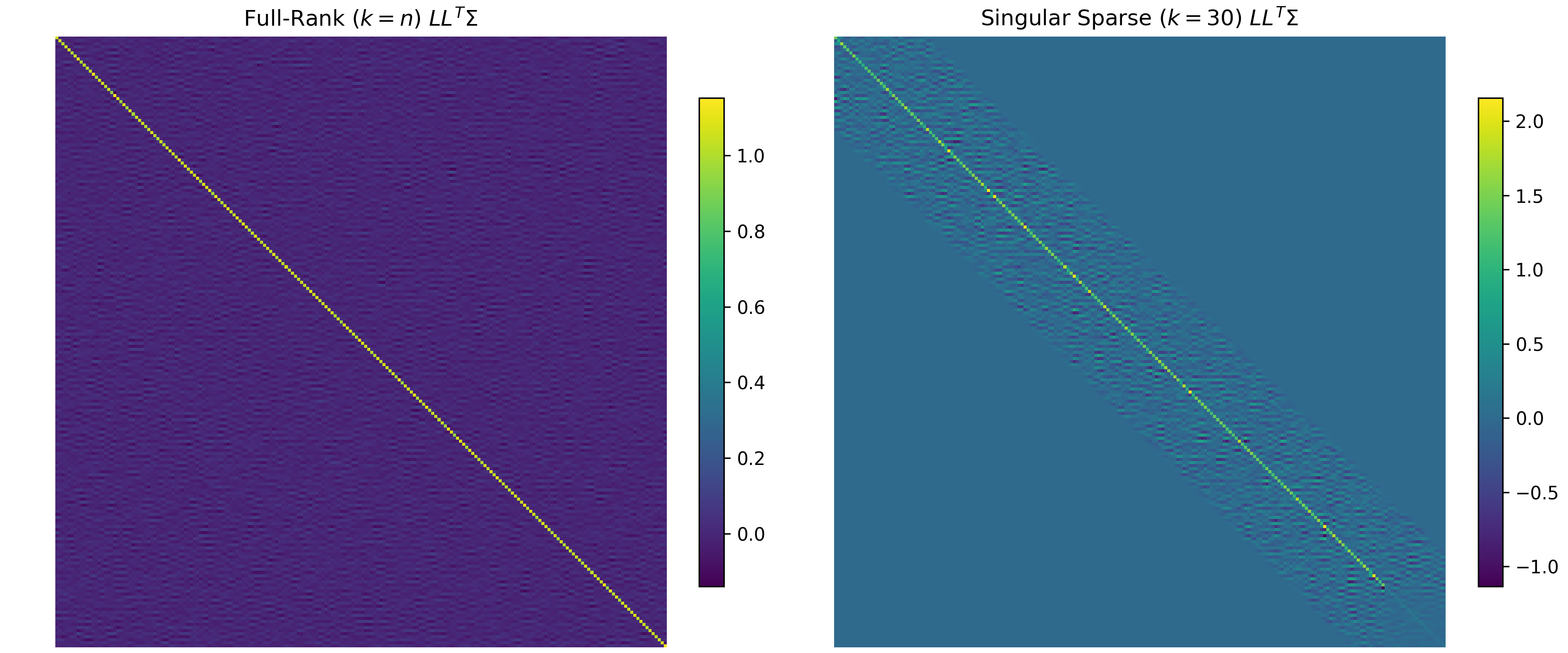}
    \caption{Full-Rank KL: 1.05 \\ Singular Case KL: 81.0}
    \label{fig:err_sparse}
\end{figure}

To compute the proxy KL divergence, we apply the unified block SLQ estimator. This computes $\frac{1}{2}\operatorname{tr}(f(L^T \Sigma L))$ by sampling subblocks of size $s \times s$ from the effective dimension of $L$, performing block Lanczos iterations, and averaging the spectral quadrature estimates. In the full-rank case, we compare this estimate against the direct KL divergence. In the singular case, only the SLQ estimate is computable. Figure~\ref{fig:err} shows the KL divergence estimates under the least-squares Cholesky solution for full-rank and singular settings. Figure~\ref{fig:err_sparse} shows the same comparison when using a sparse Cholesky solver, where the singular-case KL is significantly reduced.

\paragraph{MNIST Classification with SLQ-KL Regularization.}
We train a simple multilayer perceptron (MLP) on MNIST with an architecture consisting of a fully connected input layer (784 to 8), a ReLU activation, and a final linear layer projecting to 10 output classes. The model is trained on 10\% of the MNIST dataset and evaluated on the remaining 90\%.

We use a mini-batch size of 4 and apply the SLQ-KL penalty to the hidden activations of the first layer. At each training step, we compute the hidden activation matrix \(H \in \mathbb{R}^{4 \times 8}\), center it to obtain \(H_c\), and compute the empirical covariance matrix
\[
\Sigma = \frac{1}{3}\,H_c^T H_c \;\in\; \mathbb{R}^{8 \times 8}.
\]
To estimate the KL penalty \(\frac{1}{2} \tr(f(\Sigma))\), we apply the unified block SLQ estimator with
\[
s = \min\bigl(\mathrm{rank}(\Sigma),\,\text{hidden\_dim}\bigr), \quad t = 3,
\]
subblocks, using \(q = 2\) orthonormal block probes and \(k = 8\) Lanczos steps per subblock. The identity matrix \(L = I\) is used as the Cholesky factor. This estimate is added to the cross-entropy loss with weight \(\beta = 0.01\). The network is optimized using Adam with learning rate \(10^{-3}\) for 15 epochs. We compare the regularized model against a baseline trained with \(\beta = 0\), holding all other parameters fixed.

Figures~\ref{fig:accuracy_sql} and~\ref{fig:labels_pred} report the training curves and final test-time predictions. KL-regularized training improves generalization and corrects several classification errors made by the unregularized model.

\section{Supplementary Figures} \label{appendix:figures}

\newpage
\section{Algorithms} \label{appendix:unified_block_slq}

\subsection{Stochastic Lanczos Quadrature (SLQ)}

\begin{algorithm}[h]
\caption{Stochastic Lanczos Quadrature (SLQ) Estimator for the KL Divergence}
\label{algo}
\begin{algorithmic}[1]
\STATE \textbf{Input:} Covariance operator \(\Sigma\), Cholesky factor \(L\), number of probes \(q\), Lanczos steps \(k\).
\STATE Define \(\mathcal{A}(v) = L^T[\Sigma(Lv)]\).
\FOR{\(i = 1, \dots, q\)}
    \STATE Generate \(z_i \in \{-1,1\}^n\), normalize \(v_i = z_i/\|z_i\|\).
    \STATE Run Lanczos with \(\mathcal{A}\) starting from \(v_i\), obtain tridiagonal \(T_i\).
    \STATE Diagonalize \(T_i = U_i \operatorname{diag}(\mu_{i1}, \dots, \mu_{ik}) U_i^T\).
    \STATE Set weights \(w_{ij} = (U_i(1,j))^2\) and compute \(\eta_i = \sum_{j=1}^{k} w_{ij} f(\mu_{ij})\), where \(f(\lambda) = \lambda - \ln\lambda - 1\).
\ENDFOR
\STATE Estimate the trace: \(\widehat{\operatorname{tr}}(f(\mathcal{A})) = \frac{n}{q} \sum_{i=1}^{q} \eta_i\).
\STATE \textbf{Output:} Estimated KL divergence \(\widehat{D}_{\mathrm{KL}} = \frac{1}{2} \widehat{\operatorname{tr}}(f(\mathcal{A}))\).
\end{algorithmic}
\end{algorithm}

\vspace{0.1cm}
\begin{flushright}
	\scriptsize \framebox{\parbox{2.5in}{Government License: The
			submitted manuscript has been created by UChicago Argonne,
			LLC, Operator of Argonne National Laboratory (``Argonne").
			Argonne, a U.S. Department of Energy Office of Science
			laboratory, is operated under Contract
			No. DE-AC02-06CH11357.  The U.S. Government retains for
			itself, and others acting on its behalf, a paid-up
			nonexclusive, irrevocable worldwide license in said
			article to reproduce, prepare derivative works, distribute
			copies to the public, and perform publicly and display
			publicly, by or on behalf of the Government. The Department of Energy will provide public access to these results of federally sponsored research in accordance with the DOE Public Access Plan. http://energy.gov/downloads/doe-public-access-plan. }}
	\normalsize
\end{flushright}

\end{document}